\documentclass[10pt,a4paper]{article}
\usepackage[utf8]{inputenc}
\usepackage[T1]{fontenc}
\usepackage{lmodern} 
\usepackage[francais,english]{babel}
\usepackage{amsmath,amsthm}
\usepackage{mathrsfs}
\usepackage{bbm} 
\usepackage{amssymb}
\usepackage{amsfonts}
\usepackage{pifont}
\usepackage{stmaryrd}
\usepackage{dsfont}
\usepackage{graphicx,pstricks}
\usepackage{boites,boites_exemples}
\usepackage{enumerate}
\usepackage{lscape}
\usepackage[all]{xy}
\usepackage{hyperref}
\usepackage{fullpage}
\input xy
\xyoption{all}

\newcommand{\be}{\begin{enumerate}}
\newcommand{\ee}{\end{enumerate}}
\newcommand{\R}{\mathbb{R}}

\newcommand{\N}{\mathbb{N}}
\newcommand{\C}{\textnormal{C}}

\newcommand{\W}{\textnormal{W}}

\newcommand{\D}{\textnormal{D}}
\newcommand{\U}{\textnormal{U}}
\newcommand{\E}{\textnormal{E}}

\renewcommand{\H}{\textnormal{H}}
\renewcommand{\L}{\textnormal{L}}
\newcommand{\Ll}{\textnormal{L}_{\textnormal{loc}}}

\newcommand{\T}{\mathbb{T}}
\newcommand{\X}{\textnormal{X}}

\newcommand{\V}{\textnormal{V}}

\renewcommand{\div}{\textnormal{div}}

\newcommand{\dd}{\mathrm{d}}
\newcommand{\na}{\nabla}
\newcommand{\pa}{\partial}
\newcommand{\Ld}{\textnormal{L}}
\newcommand{\Ldiv}{\Ld_{\textnormal{div}}}

\newtheorem{thm}{Theorem}[section]

\newtheorem{coro}[thm]{Corollary}
\newtheorem{propo}[thm]{Proposition}
\newtheorem{lem}[thm]{Lemma}
\newtheorem{rem}{Remark}[section]
\newtheorem{nota}{Notation}[section]
\newtheorem{defi}{Definition}[section]

\newcommand{\eps}{\varepsilon}

\numberwithin{equation}{section}

\title{Large time behavior of small data solutions to the Vlasov-Navier-Stokes system on the whole space}
\author{Daniel Han-Kwan\footnote{Centre de Math\'ematiques Laurent Schwartz (UMR 7640), CNRS, Ecole Polytechnique, Institut Polytechnique de Paris, 91128 Palaiseau Cedex, France. (\href{mailto:daniel.han-kwan@polytechnique.edu}{daniel.han-kwan@polytechnique.edu})}}
\begin{document}
\maketitle

\begin{abstract}
We study the large time behavior of small data solutions to the Vlasov-Navier-Stokes system on $\R^3 \times \R^3$. We prove that the kinetic distribution function concentrates in velocity to a Dirac mass supported at $0$, while the fluid velocity homogenizes to $0$, both at a polynomial rate. The proof is based on two steps, following the general strategy laid out in \cite{HKMM}: (1) the energy of the system decays with polynomial rate, assuming a uniform control of the kinetic density, (2) a bootstrap argument allows to obtain such a control.
This last step requires a fine understanding of the structure of the so-called Brinkman force, which follows from a family of new identities for the dissipation (and higher versions of it) associated to the Vlasov-Navier-Stokes system.
\end{abstract}

\section{Introduction}

Consider the Vlasov-Navier-Stokes system set in $\R^3 \times \R^3$:
\begin{align}
\label{eq:vlasov}\partial_t f + v \cdot \nabla_x  f + \div_v [f(u-v)] &=0,\\
\label{eq:ns}\partial_t u + u\cdot \nabla u -  \Delta u + \nabla p &= j_f-\rho_f u, \\
\label{eq:ns2}\div \, u & = 0,
\end{align}
where 
\begin{align*}
\rho_f(t,x) &:= \int_{\R^3} f(t,x,v)\,\dd v, \\
j_f(t,x) &:= \int_{\R^3} vf(t,x,v)\,\dd v.
\end{align*}
This system aims at describing the dynamics of an aerosol, that is, loosely speaking, a cloud of fine particles immersed in a  (homogeneous, incompressible)  fluid (e.g. the air); the kinetic distribution function $f(t,\cdot,\cdot)$ describes the density  of the particles in phase space $\R^3 \times \R^3$, while the fluid is described by its velocity field $u(t,\cdot)$ and pressure scalar field $p(t,\cdot)$.
The forcing term in the Navier-Stokes equations, which accounts for the exchange of momentum between the particles and the fluid is referred to as the \emph{Brinkman force}. Several variants of the model are possible (to account for more complex physics) but the Vlasov-Navier-Stokes system stands as an important prototype to build on. See e.g. \cite{BGLM} or the introduction of \cite{HKMM} (and references therein) for more details on modelling issues and on the mathematical context. 

We define the \emph{energy} and the \emph{dissipation}  of the Vlasov-Navier-Stokes system as
\begin{align}
\label{eq:energy}
&\E(t) := \frac12\int_{\R^3}|u(t,x)|^2\,\dd x +\frac12\int_{\R^3\times\R^3} f(t,x,v)|v|^2\,\dd v\,\dd x, \\
\label{eq:dissip}
&\D(t) := \int_{\R^3\times\R^3} f(t,x,v)|u(t,x)-v|^2\,\dd v\,\dd x +  \int_{\R^3} |\nabla u(t,x)|^2 \, \dd x .
\end{align}
Formally, the following energy--dissipation identity holds:
\begin{equation}
\label{eq:introenergydissip}
\frac{\dd}{\dd t} \E(t) + \D(t) =0.
\end{equation}
We consider global weak solutions to the Vlasov-Navier-Stokes system  that satisfy an energy-dissipation inequality as built in~\cite{bou-des-grand-mou}\footnote{As a matter of fact, \cite{bou-des-grand-mou} builds such a solution on $\T^3 \times \R^3$, but the proof can be adapted to $\R^3 \times \R^3$, following the arguments explained in \cite[Appendix A]{HM3}.} (see also \cite{bougramou,HM3} for more recent developments). 
Let us recall precisely this notion.

\begin{defi}\label{def:adm} We shall say that $(f_0, u_0)$ is an admissible initial condition if 
\begin{align}
\label{eq:initu}u_0 &\in\Ldiv^2(\R^3)= \{ \U \in \Ld^2(\R^3), \, \div \, \U=0\}, \\
\label{eq:initf}0\leq f_0 &\in\Ld^1 \cap\Ld^\infty(\R^3\times\R^3),\\
\label{eq:initfbis}(x,v)\mapsto f_0(x,v)|v|^2 &\in\Ld^1(\R^3\times\R^3),\\
\int_{\R^3 \times \R^3} f_0 \, \dd v \, \dd x &=1.
\end{align}
\end{defi}

\begin{defi}\label{def:sol}
Consider an admissible initial data $(u_0,f_0)$ in the sense of Definition \ref{def:adm}. A global weak solution of the Vlasov-Navier-Stokes system with initial condition $(u_0,f_0)$ is a pair $(u,f)$ with the regularity 
\begin{align*}
u &\in \Ll^\infty(\R_+;\Ld^2(\R^3))\cap\Ll^2(\R_+;\H^1_\div(\R^3)),\\
0\leq f&\in  \Ll^\infty(\R_+;\Ld^1\cap\Ld^\infty(\R^3\times\R^3)),\\
j_f-\rho_f u &\in\Ll^2(\R_+;\H^{-1}(\R^3)), \\
\int_{\R^3 \times \R^3} f \, \dd v \, \dd x &=\int_{\R^3 \times \R^3} f_0 \, \dd v \, \dd x =1, 
\end{align*}
with $u$ being a Leray solution of \eqref{eq:ns} -- \eqref{eq:ns2} (with initial condition $u|_{t=0}=u_0$) and $f$ a renormalized solution of \eqref{eq:vlasov} (with initial condition $f|_{t=0}=f_0$), and such that the following energy inequality holds for almost all $s \geq 0$ (including $s=0$) and all $t \geq s$,
\begin{align}
\label{ineq:nrj}\E(t) + \int_s^t \D(\sigma) \, \dd \sigma \leq \E(s).
\end{align}
\end{defi}

We aim  in this paper  at describing the long time behavior of  small data solutions to the Vlasov-Navier-Stokes system. This work can be seen as another part in the series of papers \cite{GHM}, \cite{HKMM}. 
In \cite{GHM} long time behavior is studied for the system set in a 2D rectangle with partly absorbing boundary conditions. It is shown that under a geometric control condition (the so-called exit geometric condition), there exist non-trivial smooth equilibria, and these equilibria are asymptotically stable under small localized perturbations.
In \cite{HKMM}, long time behavior is studied for the system set on $\T^3 \times \R^3$ (i.e. periodic data in the space variable). Let us discuss the later in more details  in the next subsection.

\subsection{The case of $\T^3$}
In the paper \cite{HKMM}, the question of long time behavior was tackled for the system set on $\T^3 \times \R^3$. On the torus, a key object is the so-called \emph{modulated energy}, as introduced by Choi and Kwon \cite{CK}:
\begin{equation}
\mathscr{E}(t) := \frac{1}{2} \int_{\T^3 \times \R^3} f(t,x,v)|v-\langle j_f(t,x) \rangle|^2 \, \dd v\, \dd x 
+ \frac12 \int_{\T^3} |u(t,x)-\langle u(t) \rangle|^2 \, \dd x +\frac{1}{4} |\langle j_f(t) \rangle-\langle u(t) \rangle|^2,
\end{equation}
 where $\langle \cdot \rangle$ stands for the spatial mean on $\T^3$.
 
Loosely speaking, the main result of \cite{HKMM} proves that under the condition 
$$
\mathscr{E}(0) + \| u_0 \|_{\dot{\H}^{1/2}(\T^3)} \ll 1,
$$
the fluid velocity $u$ homogenizes  as $t \to +\infty $ to the constant value $\U_0 := \frac{\langle u_0 + j_{f_0} \rangle}{2}$, while the kinetic distribution function $f(t)$ concentrates in velocity to a Dirac distribution supported at $\U_0$.  Moreover the convergences are exponentially fast.

Two main ingredients are at work in the proof of this result.
\begin{enumerate}
\item Choi and Kwon proved in \cite{CK} that $\mathscr{E}$ decays exponentially fast provided that one ensures the global control $\| \rho_f \|_{\Ld^\infty(0,+\infty; \Ld^\infty( \T^3))} < +\infty$. This is based on the (formal) modulated energy--dissipation law
\begin{equation}
\frac{\dd}{\dd t} \mathscr{E}(t) + \D(t) =0.
\end{equation}
and the fact that the bound $\| \rho_f \|_{\Ld^\infty(0,+\infty; \Ld^\infty( \T^3))} < +\infty$ provides the control
$$
\mathscr{E}(t) \lesssim \D(t), \qquad \forall t \geq 0,
$$
yielding
\begin{equation}
\label{decayT3intro}
\mathscr{E}(t) \lesssim e^{-\lambda t} \mathscr{E}(0), \qquad \forall t \geq 0,
\end{equation}
for some $\lambda>0$.
\item The second ingredient is a bootstrap analysis. Thanks to a straightening change of variables in velocity inspired by Bardos and Degond \cite{BD}, the global bound $\| \rho_f \|_{\Ld^\infty(0,+\infty; \Ld^\infty( \T^3))} < +\infty$ follows from an estimate bearing on the Lipschitz semi-norm of $u$, namely 
\begin{equation}
\label{eq:smallnauintro}
\int_0^{+\infty} \| \na_x u \|_{\Ld^\infty} \, \dd t \ll 1.
\end{equation}
In \cite{HKMM}, it is shown that this control can be ensured for a class of data close to equilibrium, precisely in the sense that $\mathscr{E}(0) + \| u_0\|_{\dot{\H}^{1/2}(\T^3)} \ll 1$.

The main idea is that higher order parabolic regularity estimates for the Navier-Stokes equations (possibly with some mild polynomial growth in time) can be interpolated with the modulated energy decay estimate~\eqref{decayT3intro} to produce the estimate~\eqref{eq:smallnauintro}. Indeed, the exponential decay of $\mathscr{E}$ yields the required integrability in time, while the smallness of $\mathscr{E}(0)$ yields the required smallness.

\end{enumerate}

\subsection{Main result}

We focus in this paper on the case of $\R^3 \times \R^3$.
We shall work in a small data regime, namely we loosely speaking require that the initial kinetic distribution and fluid velocity are small in the sense that $\| u_0\|_{\dot \H^{1/2}(\R^3)} + \| f_0 \|_{\Ld^1_v(\R^3; \Ld^\infty_x(\R^3))} \ll 1$, and that the initial energy is small as well, that is $\E(0)\ll 1$. 

Let us start by recalling some notations for moments from \cite{HKMM}.
\begin{defi}\label{def:decay}
We say that an initial condition satisfies the \emph{pointwise decay assumption} of order $q>0$ if
\begin{equation*}
(x,v)\longmapsto (1+|v|^q)f_0(x,v) \in \Ld^\infty(\R^3\times \R^3),
\end{equation*}
and in that case we denote 
\begin{equation*}
N_q(f_0):=\sup_{x\in\R^3,v\in\R^3} (1+|v|^q)f_0(x,v).
\end{equation*}
\end{defi}
\begin{defi}
\label{def:moments}For all $\alpha\geq 0$ and any measurable non-negative function $\textnormal{f}: \, \R^3 \times\R^3 \to \R_+$, we set
\begin{align*}
m_\alpha \textnormal{f}(t,x) &:= \int_{\R^3} \textnormal{f}  |v|^\alpha  \, \dd v, \\
M_\alpha \textnormal{f} (t) &:= \int_{\R^3 \times \R^3} \textnormal{f}  |v|^\alpha  \, \dd v \,\dd x. 
\end{align*}
\end{defi}

The main result of this paper is stated in the following theorem.
\begin{thm}
\label{thm}
There exists $p_0>3$ such that, for all $p \in (3,p_0]$ and all $\alpha\in(0,3/2)$, there exist $\delta>0$ and an onto nondecreasing function $\Psi$ such that the following holds.
Let $(u_0, f_0)$ be an admissible initial condition satisfying
\begin{equation}
\label{eq:addi}
\begin{aligned}
&u_0 \in \H^1(\R^3) \cap B^{s_p}_{p,p}(\R^3), \qquad s_p = 2- \frac{2}{p}, \\
&M_\alpha f_0 + N_q(f_0)<+\infty, \qquad \text{for   } \alpha>3, \, q>p+3.
\end{aligned}
\end{equation}
If 
\begin{equation}
\label{eq:introsmall}
\Psi(\| u_0\|_{\H^1(\R^3) \cap B^{s_p}_{p,p}(\R^3)} + M_\alpha f_0 + N_q(f_0) +1 )  \E(0)   \leq 1, \qquad   \| f_0 \|_{\Ld^1_v(\R^3; \Ld^\infty_x(\R^3))}\leq \delta,
\end{equation}
then there exists a continuous function $\varphi_\alpha$ cancelling at $0$, such that the  global weak solution $(u,f)$ with initial condition $(u_0,f_0)$ satifies
\begin{equation}
\label{eq:introdecay}
\E(t) \leq \frac{\varphi_\alpha(\E(0))}{(1+t)^{\alpha}}, \qquad \forall t \geq 0.
\end{equation}

\end{thm}

\begin{rem}
The uniqueness of the global weak solution follows from \cite{HMMM}.
\end{rem}

\begin{rem}Note that as opposed to the torus case \cite{HKMM}, a supplementary smallness condition on the initial kinetic distribution function is required.
\end{rem}

\begin{rem} It is likely that as in \cite{HKMM}, by relying on some instantaneous parabolic smoothing mechanism for the Navier-Stokes equations, the higher regularity assumption on $u_0$ in~\eqref{eq:addi} can be  partly dispensed with; note that this would nevertheless at least still require $u_0 \in \dot{\H}^{1/2}(\R^3)$ (with small norm). We have made the choice to not dwell on this possible development, as we think it is not essential.
\end{rem}

\begin{rem}
The decay in~\eqref{eq:introdecay} does not imply that $\sqrt\E \in \L^1(0,+\infty)$, which means that this does not enter the abstract framework of \cite[Theorem 1]{Jab}.
\end{rem}

Acording to the next lemma (see e.g. \cite[Lemma 1.1]{HKMM}), the energy $\E(t)$ allows to control the Wasserstein distance $\W_1$ of $f$ to the Dirac mass in velocity supported at $0$ with density $\rho_f(t)$.
\begin{lem}
\label{lem:W1}
 For all $t \geq 0$,
\begin{equation}
\begin{aligned}
\W_1 \left(f(t), \rho_f(t) \otimes \delta_{v=0}\right) + \left\| u(t) \right\|_{\Ld^2(\R^3)}   &\lesssim ({\E}(t))^{1/2}.
\end{aligned}
\end{equation}
\end{lem}

We therefore deduce
\begin{coro}
With the same assumptions and notations as in Theorem~\ref{thm}, for all $\alpha \in (0,3/4)$, for all $t \geq 0$,
\begin{equation}
\W_1 \left(f(t), \rho_f(t) \otimes \delta_{v=0}\right) + \left\| u(t) \right\|_{\Ld^2(\R^3)}   \leq \frac{\sqrt{\varphi_{2\alpha}(\E(0))}}{(1+t)^{\alpha}}.
\end{equation}

\end{coro}

In other words, this result proves that the kinetic distribution function concentrates in velocity to a Dirac mass supported at $0$, while the fluid velocity homogenizes to $0$. 
In particular this entails that the trivial solution $(0,0)$ is Lyapunov \emph{unstable}. This is in sharp contrast with the case of other Vlasov type equations such as the Vlasov-Poisson (see e.g. \cite{BD}) or Vlasov-Maxwell (see e.g. \cite{GS}) systems.

\begin{rem}
By weak compactness, there exist a sequence of times $(t_n)$ going to infinity and an asymptotic profile $\rho^\infty(x)$ as the weak limit of $\rho_f(t_n)$ as $n \to +\infty$.
However, because of the slow polynomial decay obtained in~\eqref{eq:introdecay}, we cannot apply \cite[Proposition 3.5]{HKMM} which  would prove the uniqueness of the asymptotic profile and the convergence without requiring to take a subsequence. 
\end{rem}

We will follow the strategy  outlined in the study of the torus case. However several important differences appear.

The first step of the proof will be to obtain the conditional large time decay of the energy, which is the analogue of the aforementioned result of Choi and Kwon \cite{CK} for the torus case.
However, in $\R^3$, in the absence of a Poincar\'e inequality (for the Lebesgue measure), we cannot expect exponential decay. At best, we can hope for a polynomial decay similar to that obtained for solutions to the Stokes (or heat) equation. We will show that we can indeed almost reach such an optimal rate, despite of the presence of a forcing (the Brinkman force) in the Navier-Stokes equations.
   
To this end, we will adapt Wiegner's method \cite{Wie}  for proving large time decay for the Navier-Stokes equation with source, but with a specific analysis of the influence of the source in the precise context of  the Vlasov-Navier-Stokes system (indeed the forcing is far from  decaying fast enough to apply directly the abstract results of \cite{Wie}). This takes into account the fine structure of the system. Loosely speaking, we will take advantage of the tight links between the Brinkman force and the dissipation $\D(t)$. As in the torus case, (polynomial) decay is achieved up to an \emph{a priori} control on the moment $\rho_f$. As a byproduct of this analysis, we obtain that the dissipation somehow decays faster than the energy itself (roughly speaking, a factor $1/t$ is gained).
The latest observation will serve as a guiding line for the upcoming analysis.

The second step is the bootstrap analysis, allowing to obtain  the required control on $\rho_f$. As in the torus case, a change of variables in velocity allows to reduce the problem to proving
$$
\int_0^{+\infty} \| \na_x u \|_{\Ld^\infty(\R^3)} \, \dd t \ll 1.
$$
As the decay of the energy is only polynomial, we cannot hope to give the exact same argument as in the torus case, where exponential decay of the (modulated) energy, after interpolation, virtually allows to provide any integrability in time. This interpolation will however still be useful to obtain smallness.

The idea is as follows. Recall that for the heat equation, it is  well-known that derivatives of the solution enjoy a better decay in time than the solution itself. With this perspective in mind, we shall also prove that better decay estimates hold for derivatives in space of the solution to Navier-Stokes, despite the forcing. This requires a fine understanding of the structure of the Brinkman force,  in relation with a notion of dissipation. This will lead to a family of new identities that account for the better integrability of the dissipation and of \emph{higher} order versions of it. We expect these identities to prove useful as well in other contexts.

\bigskip

The paper is organized as follows. In Section~\ref{sec2}, we show the polynomial decay of the energy, up to a conditional bound on the density $\rho_f$. The following is dedicated to the proof of this bound with the assumptions of Theorem~\ref{thm}.
Section~\ref{sec3} provides preliminaries (mostly directly taken from~\cite{HKMM}) for a bootstrap analysis. In Section~\ref{sec4}, the aforementioned key identities explaining higher decay of higher dissipations are provided, which finally allow to carry out the bootstrap argument in Section~\ref{sec5}.

\section{Conditional large time behavior on the whole space}
\label{sec2}

The goal of this section is to show the following conditional result.

\begin{thm}
\label{thm-cond}
Let $T>0$ and assume that $\|\rho_f\|_{\Ld^\infty(0,T; \Ld^\infty(\R^3))} <+\infty$.
Then for all $\alpha \in (0,3/2)$, there exists  a continuous function $\varphi$  cancelling at  $0$ depending on $ \|\rho_f\|_{\Ld^\infty(0,T; \Ld^\infty(\R^3))} $ but independent of $T$ such that
\begin{equation}
\label{eq:decaythmcond}
\E(t) \leq \frac{\varphi(\E(0))}{(1+t)^{\alpha}}.
\end{equation}
\end{thm}
We therefore obtain, up to the control of $\rho_f$, almost the same decay as for the Navier-Stokes  without source, that is loosely speaking the same as that of the heat equation on $\R^3$ (see \cite{Wie}).

\begin{rem}
In the following of the paper, abusing notations, $\varphi$ will always stand for a function satisfying the same properties as in the statement of Theorem~\ref{thm-cond}, but may change from line to line.
\end{rem}

We shall rely on the Fourier-splitting method of Schonbeck \cite{Scho,Scho2}, developed by Wiegner \cite{Wie} and Schonbeck and Wiegner \cite{SW}. Note however that we cannot apply directly their abstract results bearing on Navier-Stokes with a source, since this would require a strong decay on this source that we cannot expect to ensure.
As already mentioned in the introduction, we will rather rely on the fine algebraic structure of the full Vlasov-Navier-Stokes system.

The Fourier-splitting method is a way to control from below the fluid dissipation by the fluid energy, modulo several corrections, using a well-chosen (time dependent) splitting of the Fourier space.

\begin{proof}

Following \cite{Wie}, given a time-dependent cut-off function $g(t)$, by Plancherel\footnote{Throughout the paper, we use the normalized version of the Fourier transform such that $\| \hat{u} \|_{\L^2(\R^3)}  = \| u \|_{\L^2(\R^3)}$ for all $u \in \L^2(\R^3)$.}, we can write
$$
\begin{aligned}
\int_{\R^3} | \na_x u |^2 \, \dd x &= \int_{\R^3} | \xi|^2 | \widehat{u} |^2 \, \dd \xi \\
&\geq  \int_{|\xi|\geq g(t)} | \xi|^2 | \widehat{u} |^2 \, \dd \xi  \\
&\geq g^2(t) \| u\|_{\Ld^2(\R^3)}^2 - g^2(t) \int_{|\xi|\leq g(t)} | \widehat{u} |^2 \, \dd \xi.
\end{aligned}
$$
On the other hand, we have
$$
\begin{aligned}
 \int_{\R^3\times \R^3} f |v-u|^2 \, \dd v \, \dd x &\geq \frac{1}{2}  \int_{\R^3\times \R^3} f |v|^2 \, \dd v \, \dd x   -   \int_{\R^3} \rho_f |u|^2 \,  \dd x  \\
&\geq \frac{1}{2}    \int_{\R^3\times \R^3} f |v|^2 \, \dd v \, \dd x   -  {\|\rho_f\|_{\Ld^\infty(0,T; \Ld^\infty(\R^3))}}  \| u\|_{\Ld^2(\R^3)}^2.
\end{aligned}
$$
Choose now $\C_0>0$ large enough so that 
\begin{equation}
\label{largeC01}
\frac{\|\rho_f\|_{\Ld^\infty(0,T; \Ld^\infty(\R^3))}}{1+ \C_0}\leq 1/2.
\end{equation}
We will also ensure that for all $t \in [0,T]$,
\begin{equation}
\label{largeC02}
\frac{g^2(t)}{1+ \C_0}\leq 1/2.
\end{equation}
A part of  the fluid-kinetic dissipation term is then used in the following way:
\begin{equation*}
\frac{g^2(t)}{1+ \C_0} \int_{\R^3\times \R^3} f |v-u|^2 \, \dd v \, \dd x 
\geq \frac{1}{2} \frac{g^2(t)}{1+ \C_0}   \int_{\R^3\times \R^3} f |v|^2 \, \dd v \, \dd x   -  g^2(t)\frac{\|\rho_f\|_{\Ld^\infty(0,T; \Ld^\infty(\R^3))}}{1+ \C_0}  \| u\|_{\Ld^2(\R^3)}^2.
\end{equation*}
We deduce the following bound from below for the dissipation:
\begin{equation*}
\begin{aligned}
\int_{\R^3} &| \na_x u |^2 \, \dd x  +   \int_{\R^3\times \R^3} f |v-u|^2 \, \dd v \, \dd x \\
&\geq \left( 1- \frac{g^2(t)}{1+ \C_0}\right)   \int_{\R^3\times \R^3} f |v-u|^2 \\
&\quad +  \frac{1}{2} \frac{g^2(t)}{1+ \C_0}   \int_{\R^3\times \R^3} f |v|^2 \, \dd v \, \dd x +  g^2(t) \left[ 1 - \frac{\|\rho_f\|_{\Ld^\infty(0,T; \Ld^\infty(\R^3))}}{1+ \C_0} \right]  \| u\|_{\Ld^2(\R^3)}^2 \\
&\quad - g^2(t) \int_{|\xi|\leq g(t)} | \widehat{u} |^2 \, \dd \xi.
\end{aligned}
\end{equation*}
Thanks to~\eqref{largeC01}--\eqref{largeC02}, we get
\begin{equation*}
\begin{aligned}
\int_{\R^3} | \na_x u |^2 \, \dd x  +   \int_{\R^3\times \R^3} f |v-u|^2 \, \dd v \, \dd x &\geq\frac{1}{2} \int_{\R^3\times \R^3} f |v-u|^2  +  \frac{1}{2} \frac{g^2(t)}{1+ \C_0}  \left[  \int_{\R^3\times \R^3} f |v|^2 \, \dd v \, \dd x +   \| u\|_{\Ld^2(\R^3)}^2\right] \\
&\quad - g^2(t) \int_{|\xi|\leq g(t)} | \widehat{u} |^2  \, \dd \xi.
\end{aligned}
\end{equation*}
We set $\tilde{g}^2(t):= \frac{1}{2} \frac{g^2(t)}{1+ \C_0} $. 
By the energy inequality, we end up with the following key inequality: for almost all $s \geq 0$ and all $s \leq t \leq  T$,
\begin{equation}
\label{key}
\begin{aligned}
\E(t) + \int_s^t \tilde{g}^2(\tau) \E(\tau)\, \dd \tau &+ \frac{1}{2} \int_s^t  \int_{\R^3\times \R^3} f |v-u|^2 \, \dd v \, \dd x \dd \tau  \\
&\leq \E(s) +  \int_s^t g^2(\tau) \int_{|\xi|\leq g(\tau)}  |\widehat{u} |^2 \, \dd \xi \dd \tau.
\end{aligned}
\end{equation}
We need to control the last term of the rhs of~\eqref{key}. To this end, as in  \cite{Wie} we use the fact that $u$ solves the Navier-Stokes equation with a source term.
Let  $U_0(t,x)$ be the solution to the heat equation in $\R^3$ starting from $u_0$ at $t=0$, i.e.
$$
\pa_t U_0 - \Delta U_0 = 0, \qquad U_0|_{t=0}= u_0.
$$
 (For later use, recall that $U_0$ decays in $\Ld^2(\R^3)$ like $1/t^{3/4}$.) Taking the Fourier transform in~\eqref{eq:ns}, we obtain
$$
\pa_t \widehat{u} + |\xi|^2 \widehat{u}  =- \widehat{u\cdot \na u}  +  \widehat{F}  + \widehat{\na P}, \qquad \widehat{u}|_{t=0}= \widehat{u_0}
$$ 
where $F= j_f - \rho_f u$. By Duhamel formula, this yields
\begin{equation}
\label{ns-fourier}
\widehat{u}= \widehat{U_0} + \int_0^t \left( - \widehat{u\cdot \na u}  +  \widehat{F}  \right) e^{(s-t)|\xi|^2} \, \dd s + 
\int_0^t \widehat{\na P} e^{(s-t)|\xi|^2} \, \dd s.
\end{equation}
Thanks to the incompressibility of $u$, we must have
$$
\widehat{\na P} =   \frac{\xi \cdot \left( - \widehat{u\cdot \na u}  +  \widehat{F}  \right)}{|\xi|^2} \xi,
$$
so that for all $\xi \in \R^3$,
$$
|\widehat{\na P} (\xi)| \leq |(- \widehat{u\cdot \na u}  +  \widehat{F})(\xi) |.
$$
Integrating~\eqref{ns-fourier} with respect to $\xi$ on the ball $\{ |\xi| \leq g(\tau)\}$, 
the outcome is
\begin{align*}
\int_{|\xi|\leq g(\tau)}  |\widehat{u}(\tau) |^2 \, \dd \xi 
&\lesssim   \| U_0(\tau) \|_{\Ld^2(\R^3)}^2 + \int_{|\xi|\leq g(\tau)}  \left(\int_0^t   \left|  \widehat{u\cdot \na u} \right| +  \left| \widehat{F}  \right| \, \dd \tau \right)^2 \dd \xi.
\end{align*}
Thanks to the incompressibility of $u$, we recall that we can write ${u\cdot \na u}  =  \div (u \otimes u)$, so that
\begin{align*}
 \int_{|\xi|\leq g(\tau)}   \left(\int_0^t   \left|  \widehat{u\cdot \na u} \right| \, \dd \tau  \right)^2 \dd \xi &\lesssim  \int_{|\xi|\leq g(\tau)}  |\xi|^2  \left(\int_0^t   \left| \widehat{u \otimes u}  \right| (\tau,\xi) \, \dd \tau  \right)^2 \dd \xi \\
 &\lesssim  g(\tau)^{2+3} \left(\int_0^t   \left\| {u \otimes u}  \right\|_{\Ld^1(\R^3)} (\tau) \, \dd \tau  \right)^2. 
\end{align*}
Likewise, we obtain
$$
 \int_{|\xi|\leq g(\tau)}   \left(\int_0^t   \left|  \widehat{F} (\tau,\xi)\right| \, \dd \tau  \right)^2 \dd \xi \lesssim
 g(\tau)^{3} \left(\int_0^t   \left\| F  \right\|_{\Ld^1(\R^3)} (\tau) \, \dd \tau  \right)^2.
$$
Therefore we have proved the estimate
\begin{equation}
\label{eq-wieg}
\begin{aligned}
\int_{|\xi|\leq g(\tau)}  |\widehat{u}(\tau) |^2 \, \dd \xi 
&\leq \C\Big(   \| U_0(\tau) \|_{\Ld^2(\R^3)}^2 + g^{5}(\tau)  \left( \int_0^\tau \| u(s)\|_{\Ld^2(\R^3)}^2 \, \dd s\right)^2 \\
&\quad + g^3(\tau) \left( \int_0^\tau \| (j_f - \rho_f u) (s)\|_{\Ld^1(\R^3)} \, \dd s\right)^2\Big).
\end{aligned}
\end{equation}
 The first two terms in the rhs of \eqref{eq-wieg} will be treated exactly as in Wiegner \cite[pp. 307-308]{Wie}. Only the last one is new. We write using Cauchy-Schwarz,
\begin{equation}
\label{eq-brink0}
 \int_0^\tau \| (j_f - \rho_f u) (s)\|_{\Ld^1(\R^3)} \, \dd s \leq \| \rho_f\|_{\Ld^\infty(0,T; \Ld^\infty(\R^3))}  \int_0^\tau \left(\int_{\R^3\times \R^3} f |v-u|^2 \, \dd v \, \dd x \right)^{1/2} \, \dd s,
\end{equation}  
where we have used the normalization
$$
\int_{\R^3 \times \R^3} f (t) \, \dd v \, \dd x =\int_{\R^3 \times \R^3} f_0 \, \dd v \, \dd x  =1.
$$
Now, we can use~\eqref{key},~\eqref{eq-wieg} and~\eqref{eq-brink0} with a Gronwall-like argument (see \cite{Wie}) that is summarized in the following statement.
  \begin{lem}
  \label{lem-gron}
  Let $y(t)$ satisfy the following differential inequality. For almost all $s \geq 0$ and all  $s \leq t \leq  T$, 
  $$
  y(t) + \int_s^t \tilde{g}^2 (\tau) y(\tau) \, \dd \tau \leq y(s) + \int_s^t \beta(\tau) \, \dd \tau.
  $$
 Then for almost all $t \in [0,T]$,
 $$
 y(t) \leq y(0) \exp\left( - \int_0^t \tilde{g}^2 (\tau) \, \dd \tau\right) + \int_0^t \exp\left( - \int_\tau^t \tilde{g}^2 (r) \, \dd r\right) \beta(\tau) \, \dd \tau.
 $$
  \end{lem}
  Applying Lemma~\ref{lem-gron} with 
 \begin{align*}
 y(t) &= \E(t), \\
 \beta(\tau) &=  - \frac{1}{4} \int_{\R^3\times \R^3} f |v-u|^2 \, \dd v \, \dd x  + \C g^2(\tau) \| U_0(\tau)\|_{\Ld^2(\R^3)}^2 \\
 &\quad + \C g^{7}(\tau)  \left( \int_0^\tau \| u(s)\|_{\Ld^2(\R^3)}^2 \, \dd s\right)^2 \\
 &\quad + \C g^5(\tau) \left( \int_0^\tau \| (j_f - \rho_f u) (s)\|_{\Ld^1(\R^3)} \, \dd s\right)^2,
 \end{align*}
 we end up with the key inequality
  \begin{equation}
  \label{eq-k}
  \begin{aligned}
  \E(t) &\exp\left(\int_0^t \tilde{g}^2(s) \, \dd s  \right)
  +  \frac{1}{2} \int_0^t  \left(\int_{\R^3\times \R^3} f |v-u|^2 \, \dd v \, \dd x\right)  \exp\left(\int_0^\tau \tilde{g}^2(r) \, \dd r  \right) \dd \tau\\
  &\leq 
\E(0)+ \C  \int_0^t g^2(\tau) \| U_0(\tau) \|_{\Ld^2(\R^3)}^2 \exp\left(\int_0^\tau \tilde{g}^2(r) \, \dd r  \right) \dd \tau \\
  &\quad + \C  \int_0^t g^7(\tau)\left( \int_0^\tau \| u(r)\|_{\Ld^2(\R^3)}^2 \, \dd r\right)^2 \exp\left(\int_0^\tau \tilde{g}^2(r) \, \dd r  \right) \dd \tau \\
&\quad + \C  \int_0^t g^5(\tau)\left(  \int_0^\tau \left(\int_{\R^3\times \R^3} f |v-u|^2 \, \dd v \, \dd x \right)^{1/2} \, \dd r \right)^2 \exp\left(\int_0^\tau \tilde{g}^2(r) \, \dd r  \right) \dd \tau.
  \end{aligned}
  \end{equation}
  We are finally in position to choose an appropriate function $g(t)$. We pick up, still following  \cite{Wie},
  $$
  \tilde{g}^2(t) = \frac{\alpha}{10+ t} \Longleftrightarrow g^2(t) = \frac{2\alpha(1+\C_0)}{10+t}, 
  $$
  with $\alpha>0$ to be carefully determined. Unlike \cite{Wie}, we will always consider $1<\alpha<3/2$, so that~\eqref{largeC02} is satisfied. One remarks that~\eqref{largeC02} can thus indeed be ensured, picking $\C_0$ large enough.
By construction, we have 
  $$
   \exp\left(\int_0^t \tilde{g}^2(s) \, \dd s  \right) = (10+t)^\alpha.
   $$
   Note that we have, as $\alpha<3/2$, 
\begin{equation}
\label{g2}
   \begin{aligned}
   \E(0)+  \int_0^t g^2(\tau) \| U_0(\tau) \|_{\Ld^2(\R^3)}^2 \exp\left(\int_0^\tau \tilde{g}^2(r) \, \dd r  \right) \dd \tau 
   &\lesssim    \E(0)+  \int_0^t \frac{\dd \tau}{(1+ \tau)^{1+ 3/2 -\alpha}} \\
   &\lesssim 1.
  \end{aligned}
\end{equation}
  We now explain the iteration procedure that allows to obtain~\eqref{eq:decaythmcond} for values of $\alpha$ that are less than but arbitrarily close to $3/2$. Assume that on $[0,T]$,
  \begin{equation}
  \label{hypbeta-1}
  \E(t) \lesssim \frac{1}{(1+t)^\beta},
\end{equation}
  with $0\leq \beta <3/2$.
As by the energy inequality, $\E(t)$ is bounded, we will be able to start later with $\beta=0$.
     For $\beta <1$, we have
  $$
   \int_0^t g^7(\tau)\left( \int_0^\tau \| u(r)\|_{\Ld^2(\R^3)}^2 \, \dd r\right)^2 \exp\left(\int_0^\tau \tilde{g}^2(r) \, \dd r  \right) \dd \tau  \lesssim
     \int_0^t (1+\tau)^{\alpha -   2(\beta-1) - 7/2} \dd \tau. 
       $$
       This is bounded by $(1+t)^{\alpha-2\beta-1/2}$ if $\alpha -2 \beta-3/2>-1$, or directly  by 
       a constant if $\alpha -2 \beta-3/2<-1$.
       For $\beta>1$, a similar computation yields directly a bound by a constant (since $\alpha<3/2$).
To summarize, assuming~\eqref{hypbeta-1}, we have 
\begin{equation}
\label{g7}
\begin{aligned}
   \int_0^t g^7(\tau)\left( \int_0^\tau \| u(r)\|_{\Ld^2(\R^3)}^2 \, \dd r\right)^2 &\exp\left(\int_0^\tau \tilde{g}^2(r) \, \dd r  \right) \dd \tau \\
    &\lesssim   (1+t)^{\alpha-2\beta-1/2} \qquad \text{if  } \beta<1,  \quad \alpha -2 \beta-3/2>-1, \\
    &\lesssim 1  \qquad \text{if  } \beta<1,  \quad \alpha -2 \beta-3/2<-1, \\
      &\lesssim 1   \qquad \text{if  } \beta>1.
      \end{aligned}
\end{equation}
 
Let us assume as well that on $[0,T]$,

  \begin{equation}
  \label{hypbeta-2}
  \int_0^t  \left(\int_{\R^3\times \R^3} f |v-u|^2 \, \dd v \, \dd x\right) (10+\tau)^\alpha \dd \tau  \lesssim \frac{(10+t)^\alpha}{(1+t)^\beta}.
  \end{equation}
  Remark that if \eqref{hypbeta-2} holds for some $\alpha$, then it holds as well for all $\widetilde\alpha\geq \alpha$.
If $2\alpha - \beta -3/2>0$,  by Cauchy-Schwarz, there holds
  $$
  \begin{aligned}
   \int_0^t g^5(\tau) &\left(   \int_0^\tau \left(\int_{\R^3\times \R^3} f |v-u|^2 \, \dd v \, \dd x \right)^{1/2} \, \dd r \right)^2 \exp\left(\int_0^\tau \tilde{g}^2(r) \, \dd r  \right) \dd \tau \\
    &\lesssim  \int_0^t g^5(\tau)\left(   \int_0^\tau \left(\int_{\R^3\times \R^3} f |v-u|^2 \, \dd v \, \dd x \right) (10+r)^\alpha \, \dd r \right) \\
    &\qquad \qquad \qquad \times\left(   \int_0^\tau \frac{1}{(10+r)^\alpha} \, \dd r\right)  \exp\left(\int_0^\tau \tilde{g}^2(r) \, \dd r  \right) \dd \tau \\
     &\lesssim    \int_0^t  (1+\tau)^{2\alpha -\beta -5/2} \, \dd\tau \lesssim (1+t)^{2\alpha-\beta-3/2}.
     \end{aligned}
  $$
  Else, if $2\alpha - \beta -3/2< 0$, we have a bound by a constant. To summarize, assuming~\eqref{hypbeta-2} we have 
\begin{equation}
\label{g5}
\begin{aligned}
\int_0^t g^5(\tau)&\left(   \int_0^\tau \left(\int_{\R^3\times \R^3} f |v-u|^2 \, \dd v \, \dd x \right)^{1/2} \, \dd r \right)^2 \exp\left(\int_0^\tau \tilde{g}^2(r) \, \dd r  \right) \dd \tau  \\
    &\lesssim   (1+t)^{2\alpha-\beta-3/2} \qquad \text{if  } 2\alpha - \beta-3/2>0, \\
    &\lesssim 1  \qquad \text{if  }   2 \alpha - \beta-3/2<0.
      \end{aligned}
\end{equation}

\bigskip

Now we argue by induction in order to increase the admissible values of $\beta$.

Start with $\beta=0$, and take $\alpha=1$. The a priori estimates~\eqref{hypbeta-1} and~\eqref{hypbeta-2} are indeed satisfied since by the energy inequality
$$
 \E(t)  +  \int_0^t  \left(\int_{\R^3\times \R^3} f |v-u|^2 \, \dd v \, \dd x\right)  \dd \tau  \lesssim 1,
  $$
  so that
  $$
  \int_0^t  \left(\int_{\R^3\times \R^3} f |v-u|^2 (10+ \tau) \, \dd v \, \dd x\right)  \dd \tau  \lesssim (10+t).
  $$
Using~\eqref{eq-k} together with~\eqref{g2}, \eqref{g7} and~\eqref{g5}, we obtain
\begin{equation*}
{(10+t)}  \E(t) +   \int_0^t  \left(\int_{\R^3\times \R^3} f |v-u|^2 \, \dd v \, \dd x\right) (10+\tau) \dd \tau \\
\lesssim  1+  (1+t)^{1-1/2} +  (1+t)^{2-3/2},
\end{equation*}
so that
\begin{equation*}
\begin{aligned}
  \E(t) &+ \frac{1}{(10+t)}   \int_0^t  \left(\int_{\R^3\times \R^3} f |v-u|^2 \, \dd v \, \dd x\right) (10+\tau) \dd \tau \\
  &\lesssim  \frac{1}{(1+t)}  + \frac{1}{(1+t)^{1/2}}  \\
  &\lesssim  \frac{1}{(1+t)^{1/2}},
  \end{aligned}
\end{equation*}
which means that~\eqref{hypbeta-1} and~\eqref{hypbeta-2} are satisfied for $\beta=1/2$, $\alpha=1$. 
%
%

\bigskip

Now we start again with $\beta_1=1/2$ and $\alpha_2> 1$ to be fixed later. We obtain (again thanks to \eqref{eq-k} together with~\eqref{g2}, \eqref{g7} and~\eqref{g5}) that
\begin{equation*}
  \E(t) + \frac{1}{(10+t)^{\alpha_1}}   \int_0^t  \left(\int_{\R^3\times \R^3} f |v-u|^2 \, \dd v \, \dd x\right) (10+\tau)^{\alpha_2} \dd \tau 
  \lesssim  \frac{1}{(1+t)^{\alpha_2}}  +  \frac{1}{(1+t)^{\beta_1 + 3/2- \alpha_2}}.
\end{equation*}
As $\alpha_2$ can be taken arbitrarily close to $1$, this yields the decay $  \frac{1}{(1+t)^{1^-}}$. 
For all $\eps\ll 1$, we can thus  find $\alpha_2>1$ such that, denoting 
$$\beta_2 = 1- \eps$$ 
the controls \eqref{hypbeta-1} and~\eqref{hypbeta-2} are satisfied for $\alpha= \alpha_2$ and $\beta=\beta_2$.
Applying again the same procedure, we deduce that for all $\alpha_3\geq \alpha_2$,
\begin{equation*}
  \E(t) + \frac{1}{(10+t)^{\alpha_2}}   \int_0^t  \left(\int_{\R^3\times \R^3} f |v-u|^2 \, \dd v \, \dd x\right) (10+\tau)^{\alpha_3} \dd \tau 
  \lesssim  \frac{1}{(1+t)^{\alpha_3}}  +  \frac{1}{(1+t)^{\beta_2 + 3/2- \alpha_3}},
\end{equation*}
which gives, taking $\alpha_3= \frac{\beta_2 + 3/2}{2}(1+\eps)$, a decay in $  \frac{1}{(1+t)^{\beta_3}}$, for $\beta_3= \frac{\beta_2 + 3/2}{2}(1-\eps) $, and thus  \eqref{hypbeta-1} and~\eqref{hypbeta-2} are satisfied for $\beta =\beta_3$ and $\alpha=\alpha_3$.

\begin{rem}
The choice $\alpha_3= \frac{\beta_2 + 3/2}{2}$ may look better, but this yields some logarithmic factors in the estimates; this is why we made this small modification with the factor $(1+\eps)$ (which can be anyway taken arbitrarily close to $1$).
\end{rem}

\bigskip

This invites to define by induction, given $\beta_n$ for $n \geq 2$,
$$
\alpha_{n+1} =  \frac{\beta_n + 3/2}{2}(1+\eps) \quad \beta_{n+1}=  \frac{\beta_n + 3/2}{2}(1-\eps),
$$ 
and we get for all $n\geq 2$, on $[0,T]$,
\begin{equation}
\label{eq-final}
\begin{aligned}
  \E(t) + \frac{1}{(10+t)^{\alpha_n}}   \int_0^t  \left(\int_{\R^3\times \R^3} f |v-u|^2 \, \dd v \, \dd x\right) (10+\tau)^{\alpha_n} \dd \tau&\lesssim  \frac{1}{(1+t)^{\beta_{n}}}. 
\end{aligned}
\end{equation}

Clearly, $(\beta_n)_{n \geq 2}$ is an increasing and bounded  sequence.  Therefore it converges
and the limit is given by 
$$
\lim_{n \to +\infty} \beta_n = \frac{3}{2}\frac{1-\eps}{1+\eps}.
$$
We deduce that $(\alpha_n)_{n \geq 2}$  is converging as well, with
$$
\lim_{n \to +\infty} \alpha_n = \frac{3}{2}. 
$$
As $\eps>0$ is arbitrary, from~\eqref{eq-final} for $n$ large enough, this yields the claimed decay estimate~\eqref{eq:decaythmcond}.

\end{proof}

We have actually also obtained in the course of the proof (see~\eqref{eq-final}) the following higher decay estimate for the fluid-kinetic dissipation. 
\begin{lem}
\label{lem-decayD2}
Under the same assumptions and notations as Theorem~\ref{thm-cond}, on $[0,T]$, for all $\alpha \in (0, 3/2)$,
\begin{equation}
\label{eq:D2superdecay}
 \int_0^t  \left(\int_{\R^3\times \R^3} f |v-u|^2 \, \dd v \, \dd x\right) (1+\tau)^{\alpha} \dd \tau \lesssim \varphi(\E(0)).
 \end{equation}
\end{lem}
This decay in time is better than expected in the sense that the decay of the energy of Theorem~\ref{thm-cond} is not sufficient to yield~\eqref{eq:D2superdecay}: loosely speaking there seems to be a gain of $1/t$. 

We can also remark that by a small modification of the above proof, we can obtain a similar higher decay for the fluid dissipation.
\begin{lem}
\label{lem-decaydxu}
Under the same assumptions and notations as Theorem~\ref{thm-cond}, on $[0,T]$, for all $\alpha \in (0, 3/2)$,
\begin{equation}
 \int_0^t  \left(\int_{\R^3}  |\na_x u|^2 \,  \dd x\right) (1+\tau)^{\alpha} \dd \tau \lesssim \varphi(\E(0)).
 \end{equation}
\end{lem}

\begin{proof}[Sketch of proof] This follows from a small modification of~\eqref{key}:
\begin{equation}
\label{key-2}
\begin{aligned}
\E(t) + \int_s^t \overline{g}^2(\tau) \E(\tau)\, \dd \tau &+ \frac{1}{2} \int_s^t  \left(\int_{\R^3\times \R^3} f |v-u|^2 \, \dd v \, \dd x + \int_{\R^3} | \na u |^2 \, \dd x \right) \dd \tau  \\
&\leq \E(s) + \frac{1}{2} \int_s^t g^2(\tau) \int_{|\xi|\leq g(\tau)}  |\widehat{u} |^2 \, \dd \xi \dd \tau.
\end{aligned}
\end{equation}
where $\overline{g}^2(\tau) := \frac{1}{4} \frac{g^2(\tau)}{1+\C_0}$.
Then the remaining of the proof is the same.

\end{proof}

Lemmas~\ref{lem-decayD2} and~\ref{lem-decaydxu} inspire two ideas that will serve as guiding lines for the upcoming analysis:
\begin{itemize}

\item The derivatives in space of the solution to the Navier-Stokes equation with the Brinkman force should decay faster than the solution itself, as it should be for parabolic equations.

\item How come the fluid-kinetic dissipation seems to decay faster than the energy? This will lead to the key identities of the upcoming Lemma~\ref{lem-gaindecay}.

\end{itemize}

\section{Preliminaries for the bootstrap}
\label{sec3}

We shall prove by a bootstrap argument that there is $\C_0>0$ such that
$$
\|\rho_f\|_{\Ld^\infty(0,+\infty; \Ld^\infty(\R^3))} \leq \C_0. 
$$
Then Theorem~\ref{thm-cond} applies for $T=+\infty$ and implies Theorem~\ref{thm}.
As the early stages of the analysis are very similar to the torus case, the content of this section is mainly taken directly from \cite{HKMM}.
Let us start by defining the characteristics curves associated to the Vlasov equation~\eqref{eq:vlasov}, that are the solutions $(\X, \V)$  to the system
\begin{equation}
\label{eq:carac}
\begin{aligned}
\dot{\X}(s;t,x,v) &= \V(s;t,x,v),\\
\dot{\V}(s;t,x,v) &= u(s,\X(s;t,x,v))-\V(s;t,x,v),
\end{aligned}
\end{equation}
with $(\X(t;t,x,v),\V(t;t,x,v))=(x,v)$. This system can be solved using the DiPerna-Lions theory \cite{DPL}, exactly as in~\cite{HKMM}. Eventually we will see that $u \in \Ld^1(0,+\infty; \W^{1,\infty}(\R^3))$, so that the classical Cauchy-Lipschitz theory actually applies.

By the method of characteristics, we can write solutions to the Vlasov equation as
\begin{equation}
f(t,x,v) = e^{3t} f_0(\X(0;t,x,v), \V(0;t,x,v)).
\end{equation}
We deduce that
\begin{align}
\label{eq:rhof}
\rho_f(t,x) &= e^{3t} \int_{\R^3} f_0(t,\X(0;t,x,v),\V(0;t,x,v))\,\dd v, \\
\label{eq:jf}
j_f(t,x) &= e^{3t} \int_{\R^3} v f_0(t,\X(0;t,x,v),\V(0;t,x,v))\,\dd v.
\end{align}

\subsection{Change of variables in velocity and bounds on moments}
The first reduction in the analysis consists in relying on a  change of variables in velocity (inspired by \cite{BD}) which directly allows to get global bounds on moments; however, such a procedure requires a control of $\|\nabla u(s)\|_{\Ld^1(0,+\infty; \Ld^\infty(\R^3))}$, in the sense that this quantity has to be small enough.

The precise statement is provided in the following lemma, whose proof can be found in~\cite[Lemma 4.4]{HKMM}.
\begin{lem}
\label{charac}
There exists $\delta_0 \in (0,1]$ such that the following holds. For any $t\geq 0$ satisfying
\begin{align}
\int_0^t \|\nabla u(s)\|_{\Ld^\infty(\R^3)}\,\dd s \leq\delta_0,
\label{ineq:nablausmall}
\end{align}
and any $x \in \R^3$, the map
\begin{equation*}
\Gamma_{t,x} : v \mapsto \V(0;t,x,v),
\end{equation*}
is a  $\mathscr{C}^1$-diffeomorphism from $\R^3$ to itself satisfying furthermore
\begin{equation}
\forall v \in \R^3,\quad|\det \D_v \Gamma_{t,x}(v) | \geq \frac{e^{3t}}{2}.
\label{eq:straight diffeo}
\end{equation}
\end{lem}
As a consequence, we deduce
\begin{lem}
\label{lem:rho} If assumption~\eqref{ineq:nablausmall} of Lemma~\ref{charac} is satisfied, we have 
\begin{equation}
\label{cont-rho1} \| \rho_f\|_{\Ld^\infty(0,t;\Ld^\infty(\R^3))} \lesssim \| f_0\|_{\Ld^1_v(\R^3; \Ld^\infty_x(\R^3))}, 
\end{equation} 
\end{lem}

\begin{proof}
By~\eqref{eq:rhof}, the change of variables $v \mapsto \Gamma_{s,x}(v)$ and~\eqref{eq:straight diffeo}, we have for almost all $s\in [0,t]$
$$
\begin{aligned}
|\rho_f (s,x) |  &\lesssim \int_{\R^3} f_0( \X(0;s,x,\Gamma_{s,x}^{-1}(w), w) \, \dd w \\
&\lesssim \| f_0\|_{\Ld^1_v(\R^3; \Ld^\infty_x(\R^3))},
\end{aligned}
$$
which entails~\eqref{cont-rho1}.
\end{proof}

\subsection{Higher order energy estimates and strong existence times}

As we need to propagate regularity for the fluid velocity, we
shall need higher (i.e. $\dot\H^1$) energy estimates  for the Navier-Stokes equations.
The following proposition can be found in \cite[Proposition 5.3]{HKMM}.

\begin{propo}\label{propo:estvns3D}
There exists a universal constant $\C_\star>0$ such that the following holds. Assume that for some $T>0$ there holds
\begin{align}
\label{ineq:smallnessvns}\|u_0\|_{\dot \H^{1/2}(\R^3)}^2  + \C_\star \int_0^T \|F(s)\|_{\dot \H^{-1/2}(\R^3)}^2\,\dd s < \frac{1}{\C_\star^2}.
\end{align}
Then one has for  all $0\leq t \leq T$ the  estimate 
\begin{align}\label{ineq:estvns3D}
\|\nabla u(t)\|_{\Ld^2(\R^3)}^2 + \int_{0}^t \|\Delta u(s)\|_{\Ld^2(\R^3)}^2\,\dd s \leq  \psi\Big(\E(0), \| u_0\|_{\H^1(\R^3)},  \sup_{[0,t]} \|\rho_f(s)\|_{\Ld^\infty(\R^3)}\Big).
\end{align}
\end{propo}
Note first that choosing $\Psi$ appropriately,  by~\eqref{eq:introsmall} and an interpolation argument, we can ensure
\begin{equation}
\label{eq:H12small}
\|u_0\|_{\dot{\H}^{1/2}(\R^3)} \leq \frac{1}{2\C_\star^2}.
\end{equation}

%
In order to use the regularization offered by Proposition~\ref{propo:estvns3D}, we will need to ensure that the smallness condition \eqref{ineq:smallnessvns} is satisfied for all times. As in \cite{HKMM}, it is convenient to introduce the following terminology.

\begin{defi}[Strong existence times]\label{def:tadm}
A real number $T\geq 0$ will be said to be a \emph{strong existence time} whenever \eqref{ineq:smallnessvns} holds. 
\end{defi}


\subsection{Local in time estimates}

Using rough estimates, it is possible to obtain local estimates in time for moments and the velocity field. This is the purpose of this subsection.
We introduce another useful notation from \cite{HKMM}.
\begin{nota}\label{nota:lesssim}
The inequality $A\lesssim_{0}B$  means
\begin{equation*}
A \leq \psi\left(\|u_0\|_{\H^{1/2}(\R^3)\cap B^{s_p,p}_p(\R^3)}+M_\alpha f_0+N_q(f_0)+\E(0) + 1 \right)B,
\end{equation*}
where $\psi:\R_+\rightarrow\R_+$ is onto, continuous and nondecreasing, and $q>4$ and $\alpha>3$ are the exponents given in the statement of Theorem~\ref{thm}. 
\end{nota}

\begin{propo}\label{propo:local}
We have $u\in\Ll^1(\R_+;\Ld^\infty(\R^3))$ and $\rho_f,j_f \in\Ll^\infty(\R_+;\Ld^\infty(\R^3))$. Moreover there exists a continuous nondecreasing function $\eta:\R_+\rightarrow\R_+$ such that 
\begin{align}
\label{eq:ulinf3}
  \|u\|_{\Ld^1(0,t;\Ld^\infty(\R^3))}&\lesssim_0 \eta(t),\\
  \label{ineq:01}   \|\rho_f\|_{\Ld^\infty(0,t;\Ld^1 \cap \Ld^\infty(\R^3))}+ \|j_f\|_{\Ld^\infty(0,t;\Ld^{3/2} \cap \Ld^\infty(\R^3))}  &\lesssim_0 \eta(t).
\end{align}
Finally, for all strong existence times $t>0$, 
\begin{equation}
\label{nablaulocal}
\nabla u \in \Ld^1(0,t; \Ld^\infty (\R^3)).
\end{equation}
\end{propo}

We refer to the proofs of \cite[Lemma 4.3, Proposition 5.1 and Corollary 6.4]{HKMM}\footnote{We take the opportunity to mention that in the proof of \cite[Corollary 6.4]{HKMM}, the final argument to reach the Lipschitz regularity  is missing: one must perform the same analysis as in the proof of  \cite[Lemma 7.3]{HKMM} but without requiring uniform in time estimates. We thank Lucas Ertzbischoff for pointing out this inaccuracy.}, which although written for the torus case, apply \emph{mutatis mutandis} to the whole space case.

Note that even if we are able to ensure that all $t$ are strong existence times, Proposition~\ref{propo:local} is not yet sufficient to obtain~\eqref{ineq:nablausmall}, as the estimates, in particular~\eqref{nablaulocal}, are not uniform with respect to $t$, and thus certainly does not imply the required smallness condition of~\eqref{ineq:nablausmall}.

Proposition~\ref{propo:local} yields
\begin{lem}\label{lem:strong1}

For all $T>0$,
\begin{equation}
\label{eq:H12-general}
\int_0^{T} \| j_f - \rho_f u \|_{\dot{\H}^{1/2}(\R^3)}^2 \, \dd t  <+\infty.
\end{equation}
Furthermore, the smallness condition of Theorem~\ref{thm} ensures that $T=1$ is a strong existence time in the sense of Definition \ref{def:tadm}. 

\end{lem}

\begin{proof}
Let $T>0$. Note that by the Sobolev embedding and the Holder inequality, for all $s \in [0,T]$,
$$
\begin{aligned}
 \|(j_f - \rho_f u)(s) \|_{\dot{\H}^{1/2}(\R^3)}^2&\lesssim \| (j_f - \rho_f u)(s) \|_{\Ld^{3/2}(\R^3)}^2 \\
 &\lesssim \D(s) \| \rho_f \|_{\Ld^3(\R^3)},
 \end{aligned}
$$
where $\D$ is the dissipation introduced in~\eqref{eq:dissip}.
Using the energy--dissipation inequality and~\eqref{ineq:01}, that yields
$$
\sup_{[0,T]}  \| \rho_f \|_{\Ld^3(\R^3)}<+\infty,
$$
we obtain
$$
\int_0^{T} \| j_f - \rho_f u \|_{\dot{\H}^{1/2}(\R^3)}^2 \, \dd t <+\infty, 
$$
which concludes the proof of the first part of the lemma.

On the other hand, for $T=1$, we have
$$
\int_0^{1} \| j_f - \rho_f u \|_{\dot{\H}^{1/2}(\R^3)}^2 \, \dd t \lesssim_0 \int_0^1 \D(t) \, \dd t  \lesssim_0 \E(0) \leq \frac{1}{3\C^3_\star},
$$
by~\eqref{eq:introsmall}, choosing $\Psi$ appropriately. Recalling~\eqref{eq:H12small}, this yields that $t=1$ is a strong existence time.

\end{proof}

\bigskip

We are finally in position to set up the bootstrap argument. To this end, introduce
\begin{align}
\label{def:tstar}t^\star := \sup\left\{\text{strong existence times }t\text{ such that  }\, \int_0^t \|\nabla u(s)\|_{\Ld^\infty(\R^3)}\,\dd s < \delta_0\right\}. 
\end{align}
By Lemma~\ref{lem:strong1} and~\eqref{nablaulocal}, we must have $t^\star>0$.
The goal is to prove that $t^\star = +\infty$, which will allow, applying Lemma~\ref{charac} and Theorem~\ref{thm-cond}, to conclude the proof of Theorem~\ref{thm}.

By contradiction, we shall assume from now on that $t^\star<+\infty$.

\section{Higher decay of higher dissipation}
\label{sec4}

We fix $T \in (0, t^\star)$.
In particular, by definition of $t^\star$, we note that the characteristics~\eqref{eq:carac} are classically defined.

\bigskip

We shall study in the section what we call the \emph{higher} fluid-kinetic dissipation.

\begin{defi} 
 Let $p\geq2$. The higher fluid-kinetic dissipation (of order $p$) is the functional 
\begin{equation}
\D_p(t) := \int_{\R^3\times \R^3} f (t,x,v) |v-u(t,x)|^p \, \dd v \, \dd x .
\end{equation}
(Note that $p=2$ corresponds to the usual fluid-kinetic dissipation term $\D$ in~\eqref{eq:dissip}.)
\end{defi}

This quantity is useful to estimate the Brinkman force because of the following elementary estimate.
\begin{lem}
\label{lem-Dp}
Let $p>1$.
On $[0,T]$,
\begin{equation}
\|( j_f - \rho_f u )(t)\|_{\Ld^p(\R^3)} \leq   \| \rho_f \|_{\Ld^\infty(0,T;\Ld^\infty(\R^3))}^{\frac{p-1}{p}} \D_p^{1/p}(t).
\end{equation}
\end{lem}

\begin{proof}
By the Holder inequality, we have
\begin{align*}
\| j_f - \rho_f u \|_{\Ld^p(\R^3)}^p \leq  \| \rho_f \|_{\Ld^\infty(0,T;\Ld^\infty(\R^3))}^{p-1} \left(\int_{\R^3 \times \R^3} f |v-u|^p \, \dd v \, \dd x\right),
\end{align*}
and the lemma follows.
\end{proof}

\begin{rem} On the torus \cite{HKMM}, it turns out to be sufficient to use the rough bound
$$
\|( j_f - \rho_f u )(t)\|_{\Ld^p}  \leq \| j_f(t)\|_{\Ld^p} +  \| \rho_f u(t)\|_{\Ld^p}.
$$
In the whole space case, this is not sufficient to close the analysis, which explains why we need a finer understanding of the Brinkman force.
\end{rem}

Higher decay of the higher dissipation $\D_p$ comes from the following key identity. We will (towards the end of the bootstrap analysis) obtain that $\D_p$ for $p>2$ enjoys  a somewhat better decay than that of $\D_2$ in Lemma~\ref{lem-decayD2}.

\begin{lem}
\label{lem-gaindecay}
Let $\varphi \in \mathscr{C}^1([0,+\infty))$.
For all $p\geq 2$, and all $\gamma \in  \R$, and all $t \geq 0$, 
\begin{equation}
\label{eq-keyDpvarphi}
\begin{aligned}
\int_0^t \D_p (s) \varphi(s) \, \dd s &= \frac{1}{p} \int_0^t \D_p (s) \varphi'(s) \, \dd s \\
&\quad -\int_0^t \int_{\R^3\times \R^3} f(s,x,v) \left[ \pa_s u + (\na_x u) v \right] \cdot [v-u(s,x)] |v-u(s,x)|^{p-2} \varphi(s) \, \dd v \, \dd x \, \dd s\\
&\quad - \left[ \frac{\varphi(s)}{p} \int_{\R^3\times \R^3} f(s,x,v) |v-u(s,x)|^p \, \dd v \, \dd x \right]^t_0.
\end{aligned}
\end{equation}
\end{lem}

\begin{proof}
Write by the method of characteristics and a change of variables
\begin{align*}
\D_p (s) &= \int_{\R^3\times \R^3} f (s,x,v) |v-u|^p \, \dd v \, \dd x  \\
&=  e^{3s}  \int_{\R^3\times \R^3} f_0 (\X(0;s,x,v),\X(0;s,x,v)) |v-u|^p \, \dd v \, \dd x \\
&=\int_{\R^3\times \R^3} f_0 (x,v) \left|\V(s;0,x,v) - u (t, \X(s;0,x,v) ) \right|^p \, \dd v \, \dd x.
\end{align*}
Remark then that
\begin{align*}
\frac{\dd}{\dd s}& |\V(s;0,x,v)- u(s,\X(s;0,x,v))|^p \\
&= p \frac{\dd}{\dd s} [\V(s;0,x,v)- u(s,\X(s;0,x,v))] \cdot [\V(s;0,x,v) - u(s,\X(s;0,x,v))] \\
&\qquad \qquad \times |
\V(s;0,x,v) - u(s,\X(s;0,x,v))|^{p-2} \\
&= p \left[ u(s,\X(s;0,x,v)) - \V(s;0,x,v) -\frac{\dd}{\dd s} ( u(s,\X(s;0,x,v))) \right] \\
&\qquad \qquad \cdot [\V(s;0,x,v) - u(s,\X(s;0,x,v))] |\V(s;0,x,v) - u(s,\X(s;0,x,v))|^{p-2}.
\end{align*}
Consequently, we have
\begin{align*}
 \left|\V(s;0,x,v) - u (s, \X(s;0,x,v)) \right|^p  
 &=  - \frac{1}{p}\frac{\dd}{\dd s} |\V(s;0,x,v) - u(s,\X(s;0,x,v))|^p \\
 & \quad - \frac{\dd}{\dd s} ( u(s,\X(s;0,x,v))) \cdot  [\V(s;0,x,v) - u(s,\X(s;0,x,v))] \\
 &\qquad \qquad \qquad \times |\V(s;0,x,v) - u(s,\X(s;0,x,v))|^{p-2}.
\end{align*}
We deduce the claimed identity by integration by parts in time.
\end{proof}
In the following, we will apply this lemma for $\varphi(s)= (1+s)^{p\gamma}$,
which leads to
\begin{equation}
\label{eq-keyDp}
\begin{aligned}
&\int_0^t \D_p (s) (1+s)^{p \gamma} \, \dd s \\
&= \gamma \int_0^t \D_p (s) (1+s)^{p \gamma-1} \, \dd s \\
&\quad -\int_0^t \int_{\R^3\times \R^3} f(s,x,v) \left[ \pa_s u + (\na_x u) v \right] \cdot [v-u(s,x)] |v-u(s,x)|^{p-2} (1+s)^{p\gamma} \, \dd v \, \dd x \, \dd s\\
&\quad - \left[ \frac{(1+s)^{p\gamma}}{p} \int_{\R^3\times \R^3} f(s,x,v) |v-u(s,x)|^p \, \dd v \, \dd x \right]^t_0.
\end{aligned}
\end{equation}

Let us comment on the key identity~\eqref{eq-keyDp}: it shows that the higher dissipation integrated against a polynomial weight in time can be decomposed as a sum of:

\begin{itemize}
\item a term of the same form involving a lower order weight in time;

\item a second term involving $\pa_s u$ and $\na_x u$ that will be somehow absorbed (see Lemma~\ref{lem-bound2} below);

\item a non-negative term and a last one independent of time involving only the initial data.
\end{itemize}

Let us right away proceed with the estimate the second term of the rhs of~\eqref{eq-keyDp}.
\begin{lem}
\label{lem-bound2}
For all $p\geq 2$, and all $\gamma \in  \R$,
and all $t \geq 0$, 
\begin{align*}
&\left|\int_0^t \int_{\R^3 \times \R^3} f(s,x,v) \pa_s u  \cdot [v-u(s,x)] |v-u(s,x)|^{p-2} (1+s)^{p\gamma} \, \dd v \, \dd x \, \dd s\right|  \\
&\qquad\qquad \lesssim  \| \rho_f \|_{\Ld^\infty(0,t; \Ld^\infty(\R^3))}^{1/p}  \| (1+s)^\gamma \pa_s u \|_{\Ld^p(0,t; \Ld^p(\R^3))} \left(\int_0^t \D_p (s) (1+s)^{p \gamma} \, \dd s\right)^{\frac{p-1}{p}}, \\
&\left|\int_0^t \int_{\R^3 \times \R^3} f(s,x,v)  (\na_x u) v \cdot [v-u(s,x)] |v-u(s,x)|^{p-2} (1+s)^{p\gamma} \, \dd v \, \dd x \, \dd s\right|  \\
&\qquad\qquad \lesssim \| (1+s)^\gamma |\na_x u| m_p^{1/p} \|_{\Ld^p(0,t; \Ld^p(\R^3))} \left(\int_0^t \D_p (s) (1+s)^{p \gamma} \, \dd s\right)^{\frac{p-1}{p}},
\end{align*}
where we recall the notation
$$
m_p(s,x) = \int_{\R^3} f (s,x,v) |v|^p \, \dd v.
$$

\end{lem}

\begin{proof}This is a consequence of the Holder inequality; details are omitted.

\end{proof}
We deduce
\begin{coro}
\label{coro-Dpmain}
For all $p\geq 2$,  all $\gamma \in  \R$, all $k \in \N$,
and all $t \geq 0$, 
\begin{equation}
\label{higherDp}
\begin{aligned}
\int_0^t \D_p (s) (1+s)^{p \gamma} \, \dd s  &\lesssim \int_0^t \D_p (s) (1+s)^{p \gamma-k} \, \dd s \\
&\quad+ \| \rho_f \|_{\Ld^\infty(0,t; \Ld^\infty(\R^3))}  \| (1+s)^\gamma \pa_s u \|_{\Ld^p(0,t; \Ld^p(\R^3))}^p \\
&\quad+ \| (1+s)^\gamma |\na_x u| m_p^{1/p} \|_{\Ld^p(0,t; \Ld^p(\R^3))}^p \\
&\quad+ \int_{\R^3 \times \R^3} f_0 |v-u_0|^p \, \dd v \, \dd x.
\end{aligned}
\end{equation}
\end{coro}

\begin{proof}This follows from a combination of the last two lemmas.
We argue by induction. For $k=0$, the estimate~\eqref{higherDp} is tautological.
Assume~\eqref{higherDp} holds for some $k \in \N$. By Lemma~\ref{lem-gaindecay}, we have the identity
\begin{equation*}
\begin{aligned}
&\int_0^t \D_p (s) (1+s)^{p \gamma-k} \, \dd s \\
&= \gamma \int_0^t \D_p (s) (1+s)^{p \gamma-(k+1)} \, \dd s \\
&\quad -\int_0^t \int_{\R^3\times \R^3} f(s,x,v) \left[ \pa_s u + (\na_x u) v \right] \cdot [v-u(s,x)] |v-u(s,x)|^{p-2} (1+s)^{p\gamma-k} \, \dd v \, \dd x \, \dd s\\
&\quad - \left[ \frac{(1+s)^{p\gamma-k}}{p} \int_{\R^3\times \R^3} f(s,x,v) |v-u(s,x)|^p \, \dd v \, \dd x \right]^t_0.
\end{aligned}
\end{equation*}
Applying Lemma~\ref{lem-bound2}, we deduce the bound
\begin{equation*}
\begin{aligned}
\int_0^t \D_p (s) (1+s)^{p \gamma-k} \, \dd s 
&\lesssim  \int_0^t \D_p (s) (1+s)^{p \gamma-(k+1)} \, \dd s \\
&\quad +\| \rho_f \|_{\Ld^\infty(0,t; \Ld^\infty(\R^3))}^{1/p}  \| (1+s)^\gamma \pa_s u \|_{\Ld^p(0,t; \Ld^p(\R^3))} \left(\int_0^t \D_p (s) (1+s)^{p \gamma-k} \, \dd s\right)^{\frac{p-1}{p}} \\
&\quad +\| (1+s)^\gamma |\na_x u| m_p^{1/p} \|_{\Ld^p(0,t; \Ld^p(\R^3))} \left(\int_0^t \D_p (s) (1+s)^{p \gamma-k} \, \dd s\right)^{\frac{p-1}{p}} \\
&\quad + \int_{\R^3 \times \R^3} f_0 |v-u_0|^p \, \dd v \, \dd x.
\end{aligned}
\end{equation*}
By Young's inequality, we end up with 
\begin{equation*}
\begin{aligned}
\int_0^t \D_p (s) (1+s)^{p \gamma-k} \, \dd s 
&\lesssim  \int_0^t \D_p (s) (1+s)^{p \gamma-(k+1)} \, \dd s \\
&\quad +\| \rho_f \|_{\Ld^\infty(0,t; \Ld^\infty(\R^3))}  \| (1+s)^\gamma \pa_s u \|_{\Ld^p(0,t; \Ld^p(\R^3))}^p  \\
&\quad +\| (1+s)^\gamma |\na_x u| m_p^{1/p} \|_{\Ld^p(0,t; \Ld^p(\R^3))}^p 
 + \int_{\R^3 \times \R^3} f_0 |v-u_0|^p \, \dd v \, \dd x,
\end{aligned}
\end{equation*}
yielding~\eqref{higherDp} at rank $k+1$. We can therefore conclude by induction.
\end{proof}

\section{The bootstrap argument}
\label{sec5}

\subsection{Weighted $\Ld^2$ maximal parabolic regularity estimates}


We shall rely on  maximal parabolic regularity to get weighted in time estimates for $\pa_t u$ and $\Delta u$. Maximal regularity for the Stokes equation reads as follows (see e.g. \cite{Sal} that concerns the heat equation, but applies to Stokes after application of the Leray projection):
\begin{thm}
\label{thm-maximal}
Let $\U_0 \in \mathscr{S}'(\R^3)$ with $\div \, \U_0=0$. Let $\U$ solve the Stokes equation with a source $S$ and initial condition $\U_0$:
$$
\begin{aligned}
\pa_t \U - \Delta \U + \na p  &=S, \\
\div \,  \U &= 0, \\
\U|_{t=0} &=\U_0.
\end{aligned}
$$
For all $p,q \in (1,+\infty)$, there holds
\begin{equation}
\| \pa_t \U \|_{\Ld^p(0,+\infty; \Ld^q(\R^3))} + \| \Delta \U \|_{\Ld^p(0,+\infty; \Ld^q(\R^3))} \lesssim \| S\|_{\Ld^p(0,+\infty; \Ld^q(\R^3))} + \| \U_0 \|_{B^{s_p}_{q,p}(\R^3)},
\end{equation}
with
\begin{equation}
s_p = 2-\frac{2}{p}.
\end{equation}

\end{thm}
 
The general principle to get weighted in time estimates for the fluid velocity field $u$ will be to write that given some power $r>0$, defining
$$
\U(t,x) := (1+t)^r u(t,x),
$$
$\U$ satisfies the Stokes equation
\begin{equation}
\label{eq-weight}
\begin{aligned}
\pa_t \U - \Delta \U + \na p  &=(1+t)^r (j_f - \rho_f u)  - (1+t)^r u\cdot \na_x u  + r (1+t)^{r-1}  u, \\
\div \,  \U &= 0, \\
\U|_{t=0} &=u_0.
\end{aligned}
\end{equation}

Recall that we have fixed $T \in (0, t^\star)$. We first get the following regularity statement.
\begin{lem}
\label{lem-max1}
For all $\gamma \in (0,3/4)$, we have the estimate
\begin{equation}
\label{eq-max1}
\| (1+t)^{\gamma} \pa_t u \|_{\Ld^2(0,T; \Ld^2(\R^3))} + \| (1+t)^{\gamma} \Delta u \|_{\Ld^2(0,T; \Ld^2(\R^3))} \lesssim_0 \varphi(\E(0)) + \| u_0 \|_{\H^1(\R^3)}.
\end{equation}
\end{lem}

\begin{proof}
Let $\gamma \in (0,3/4)$. By Lemma~\ref{lem-Dp} with $p=2$ and Lemma~\ref{lem-decayD2} with $\alpha =2\gamma$, we have
\begin{equation}
\label{eq-brink1}
\| (1+t)^\gamma (j_f - \rho_f u) \|_{\Ld^2(0,T; \Ld^2(\R^3))}\lesssim \varphi(\E(0)).
\end{equation}
Furthermore, by Theorem~\ref{thm-cond}, we have for all $t \in [0,T]$,
\begin{equation}
\label{eq-ufirst}
(1+t)^{\gamma-1}  \| u (t)\|_{\Ld^2(\R^3)} \lesssim \frac{\varphi(\E(0))}{(1+t)},
\end{equation}
which is clearly uniformly (i.e. independently of $T$) bounded in $\Ld^2(0,T)$. Moreover, we estimate using the Holder inequality
$$
\| u\cdot \na_x u \|_{\Ld^2(\R^3)} \leq \| u\|_{\Ld^6(\R^3)} \| \na_x u \|_{\Ld^3(\R^3)} \leq \| u\|_{\Ld^6(\R^3)} \| \na_x u \|_{\Ld^2(\R^3)}^{1/2} \| \na_x u \|_{\Ld^6(\R^3)}^{1/2}.
$$ 
By  Sobolev embedding, we recall that $\| u\|_{\Ld^\infty(0,T;\Ld^6(\R^3))} \lesssim \| \na u\|_{\Ld^\infty(0,T;\Ld^2(\R^3))}$ and $\| \na u\|_{\Ld^2(0,T;\Ld^6(\R^3))} \lesssim \| \na_x^2 u\|_{\Ld^2(0,T;\Ld^2(\R^3))}$. Therefore, by Proposition~\ref{propo:estvns3D}, we get
\begin{equation}
\label{eq-try1}
\begin{aligned}
\| (1+t)^{\gamma/2} &u\cdot \na_x u \|_{\Ld^2(0,T; \Ld^2(\R^3))}  \\
&\lesssim \|  u \|_{\Ld^\infty(0,T; \Ld^6(\R^3))} \| (1+t)^{\gamma} \na_x u \|_{\Ld^2(0,T;\Ld^2(\R^3))}^{1/2}  \|  \na_x u \|_{\Ld^2(0,T;\Ld^6(\R^3))}^{1/2} \\
&\lesssim_0 \| (1+t)^{\gamma} \na_x u \|_{\Ld^2(0,T;\Ld^2(\R^3))}^{1/2}  \\
&\lesssim_0  \varphi(\E(0)),
\end{aligned}
\end{equation}
where we have used in the last line Lemma~\ref{lem-decaydxu} for $\alpha = 2 \gamma$.

We now set $\U = (1+t)^{\gamma/2} u$, so that $\U$ solves~\eqref{eq-weight} with $r= \gamma/2$.
By the maximal parabolic regularity result for $p=q=2$ of Theorem~\ref{thm-maximal}, we deduce
\begin{equation*}
\| \pa_t \U \|_{\Ld^2(0,T; \Ld^2(\R^3))} + \| \Delta \U \|_{\Ld^2(0,T; \Ld^2(\R^3))} \lesssim \| S\|_{\Ld^2(0,T; \Ld^2(\R^3))} + \| u_0 \|_{\H^1(\R^3)},
\end{equation*}
with $S= (1+t)^{\gamma/2} (j_f - \rho_f u)  - (1+t)^{\gamma/2} u\cdot \na_x u  + \frac{\gamma}{2} (1+t)^{\gamma/2-1}  u$.
Using~\eqref{eq-brink1}, \eqref{eq-ufirst} and \eqref{eq-try1}, we have the bound
\begin{equation}
\| S \|_{\Ld^2(0,T; \Ld^2(\R^3))} \lesssim   \varphi(\E(0)).
\end{equation}
As 
$$
\pa_t \U = (1+t)^{\gamma/2} \pa_t u  + \frac{\gamma}{2} (1+t)^{\gamma/2-1}  u ,
$$
gathering all pieces together we finally obtain
\begin{equation}
\label{eq-veryfirst}
\| (1+t)^{\gamma/2} \pa_t u \|_{\Ld^2(0,T; \Ld^2(\R^3))} + \| (1+t)^{\gamma/2} \Delta u \|_{\Ld^2(0,T; \Ld^2(\R^3))} 
\lesssim_0  \varphi(\E(0)) +  \| u_0 \|_{\H^1(\R^3)}.
\end{equation}
We have therefore obtained a better control on $u$ than~\eqref{ineq:estvns3D}  but we have not yet reached~\eqref{eq-max1}. The procedure has to be reiterated.
By Sobolev embedding, \eqref{eq-veryfirst} yields 
$$
 \| (1+t)^{\gamma/2} \nabla_x u \|_{\Ld^2(0,T; \Ld^6(\R^3))} \lesssim_0 \varphi(\E(0)) +  \| u_0 \|_{\H^1(\R^3)}.
$$
We use this new information as follows: arguing as in~\eqref{eq-try1}, we get
\begin{equation}
\label{eq-try2}
\begin{aligned}
\| (1+t)^{3\gamma/4} &u\cdot \na_x u \|_{\Ld^2(0,T; \Ld^2(\R^3))}  \\
&\lesssim \|  u \|_{\Ld^\infty(0,T; \Ld^6)} \| (1+t)^{\gamma} \na_x u \|_{\Ld^2(0,T;\Ld^2)}^{1/2}  \| (1+t)^{\gamma/2}  \na_x u \|_{\Ld^2(0,T;\Ld^6)}^{1/2} \\
&\lesssim_0  \varphi(\E(0)) +  \| u_0 \|_{\H^1(\R^3)}. 
\end{aligned}
\end{equation}
As a consequence, applying~\eqref{eq-try2} instead of~\eqref{eq-try1}, we obtain instead of~\eqref{eq-veryfirst} the enhanced inequality
\begin{equation}
\label{eq-veryfirst2}
\| (1+t)^{3\gamma/4} \pa_t u \|_{\Ld^2(0,T; \Ld^2(\R^3))} + \| (1+t)^{3\gamma/4} \Delta u \|_{\Ld^2(0,T; \Ld^2(\R^3))} 
\lesssim_0  \varphi(\E(0)) +  \| u_0 \|_{\H^1(\R^3)}.
\end{equation}
Applying recursively this procedure, this results in the bound, for all $n \geq 1$,
\begin{equation}
\label{eq-veryfirst2-reloaded}
\begin{aligned}
\| (1+t)^{\gamma \sum_{k=1}^{n} \frac{1}{2^k}} \pa_t u \|_{\Ld^2(0,T; \Ld^2(\R^3))} &+ \| (1+t)^{\gamma \sum_{k=1}^{n} \frac{1}{2^k} } \Delta u \|_{\Ld^2(0,T; \Ld^2(\R^3))} 
\\
&\lesssim_0  \varphi(\E(0)) +  \| u_0 \|_{\H^1(\R^3)}.
\end{aligned}
\end{equation}
As  $\displaystyle\sum_{k=1}^{+\infty} \frac{1}{2^k} = 1$, we may apply this procedure as many times as necessary to reach a power $\tilde\gamma$ arbitrarily close to $\gamma$: as $\gamma$ can itself be taken arbitrarily close to $3/4$, the claimed estimate follows.
\end{proof}

As a  consequence of Lemma~\ref{lem-max1}, we deduce the following control on $\| u\|_{\Ld^\infty(\R^3)}$:
\begin{coro}
\label{coro-para1}
For all $\gamma \in (0,3/4)$, there holds
\begin{equation}
\label{eq:coro-para1}
  \| (1+t)^{\gamma} u \|_{\Ld^{8/3}(0,T; \Ld^\infty(\R^3))} \lesssim_0  \varphi(\E(0)).
\end{equation}
\end{coro}

\begin{proof}
By the Gagliardo-Nirenberg inequality, 
$$
\| u \|_{\Ld^\infty(\R^3)} \lesssim   \| \Delta u \|_{\Ld^2(\R^3)}^{\alpha} \| u \|_{\Ld^{2}(\R^3)}^{1-\alpha},
$$
with
$$
0 = 0 + \left( \frac{1}{2} - \frac{2}{3} \right) \alpha + \frac{1-\alpha}{2} \Longleftrightarrow \alpha =  \frac{3}{4}.
$$
The estimate~\eqref{eq:coro-para1} then follows from the estimates of  Theorem~\ref{thm-cond} and Lemma~\ref{lem-max1}.
\end{proof}

A first application of Corollary~\ref{coro-para1} is a pointwise in time control of the $\Ld^\infty(\R^3)$ norm of the moment $m_p$, which improves as $p$ increases.
\begin{lem}
\label{lem-decaymp}
For all $p\geq 0$, $q >p+3$ and $\gamma \in (0,3/4)$, for all $t \in [0,T)$, there holds 
\begin{equation}
\label{eq-boundmp}
\| m_p(t) \|_{\Ld^\infty(\R^3))} \lesssim_0 \frac{N_q(f_0)}{(1+t)^{\gamma p}}.
\end{equation}
In particular, if $N_q(f_0)<+\infty$,  for all $\tilde\gamma\geq 0$, there is $k \in \N$ large enough, so that
\begin{equation}
\int_0^t m_p(s) (1+s)^{p\tilde\gamma-k} \, \dd s \lesssim_0 1. 
\end{equation}

\end{lem}

\begin{proof}
By the method of characteristics, we write
\begin{align*}
m_p(t,x) &= \int_{\R^3} f(t,x,v) |v|^p \, \dd v \\
&= e^{3t} \int_{\R^3} f_0(\X(0;t,x,v),\V(0;t,x,v)) |v|^p \, \dd v \
\end{align*}
By Lemma~\ref{charac}, we can use the change of variables $w:= \V(0;t,x,v) (= \Gamma_{t,x}(v))$, that yields
$$
|m_p(t,x)| \lesssim \int_{\R^3} f_0(\X(0;t,x,\Gamma_{t,x}^{-1}(w)),w) |\Gamma_{t,x}^{-1}(w)|^p \, \dd w.
$$
By~\eqref{eq:carac} we infer
$$
|\Gamma_{t,x}^{-1}(w)| \leq e^{-t} |w| + \int_0^t e^{\tau-t} \left|u\left(\tau,  \X\left(\tau;t,x,\Gamma_{t,x}^{-1}(w)\right)\right)\right| \, \dd\tau,
$$
and we deduce that
\begin{align*}
|\Gamma_{t,x}^{-1}(w)| \leq e^{-t} |w| +  \int_0^t e^{\tau-t} \| u\|_{\Ld^\infty(\R^3)} \, \dd \tau.
\end{align*}
Let $\alpha \in (0,3/4)$. By the Holder inequality and Corollary~\ref{coro-para1}, we have
\begin{align*}
\int_0^t e^{\tau-t} \| u\|_{\Ld^\infty(\R^3)} \, \dd \tau &\lesssim \left(\int_0^t \frac{e^{\frac{8}{5}(\tau-t)}}{(1+\tau)^{\frac{8}{5}\alpha}} \, \dd \tau\right)^{5/8}  \| (1+t)^{\alpha} u \|_{\Ld^{8/3}(0,T; \Ld^\infty(\R^3))} \\
&\lesssim_0 \frac{1}{(1+t)^{\alpha}}.
\end{align*}
Consequently,
$$
|m_p(t,x)| \lesssim_0 \frac{N_q(f_0)}{(1+t)^\alpha}  \int_{\R^3} \frac{|w|^p}{1+|w|^q}  \, \dd w,
$$
and the integral is finite pour $q>p+3$, hence proving~\eqref{eq-boundmp}. The other statement is just a matter of taking $k$ large enough to ensure integrability in time.

\end{proof}

%
%

\subsection{First $\Ld^p$ bounds on the source term}

In view of a subsequent application of Theorem~\ref{thm-maximal}, let us first prove some bounds in $\Ld^p$ for $p>3$ on the terms $(1+t)^{\gamma} u\cdot \na_x u  + {\gamma} (1+t)^{\gamma-1}  u$ in the source term of the Stokes equation.
\begin{lem} 
\label{lem-source1}
Let $p>3$. For all $\gamma \in \left(0,  \frac{17}{8} - \frac{7}{4p} \right)$, we have
\begin{equation}
\| (1+t)^{\gamma-1} u \|_{\Ld^p(0,T; \Ld^p(\R^p))} \lesssim_0  \varphi(\E(0)).
\end{equation}

\end{lem}

\begin{proof}
By interpolation, we write 
$$
\|u \|_{\Ld^p(\R^3)}^p \lesssim \| u \|_{\Ld^2(\R^3)}^{2}  \|u \|_{\Ld^\infty(\R^3)}^{p-2},
$$
so that by Theorem~\ref{thm-cond}, for all $\gamma >0$,
$$
\begin{aligned}
(1+t)^{p(\gamma-1)} \|u \|_{\Ld^p(\R^3)}^p &\lesssim \frac{ \varphi(\E(0)) }{(1+t)^{\frac{3}{2}^-}} (1+t)^{p(\gamma-1)}  \|u \|_{\Ld^\infty(\R^3)}^{p-2} \\
&\lesssim  \varphi(\E(0)) (1+t)^{p(\gamma-1)    -p\frac{3}{4}^- } [ (1+t)^{\frac{3}{4}^-} \| u \|_{\Ld^\infty(\R^3)} ]^{p-2},
\end{aligned}
$$ 
where $\frac{3}{4}^-$ stands for some $\alpha \in (0,3/4)$ arbitrarily close to $3/4$.
By the Holder inequality, we deduce that 
$$
\begin{aligned}
\| (1+t)^{\gamma-1} &u \|_{\Ld^p(0,T; \Ld^p(\R^p))}^p  \\
&\lesssim  \varphi(\E(0)) \| (1+t)^{p(\gamma-1)    -p\frac{3}{4}^- }\|_{\Ld^{\frac{8}{14-3p}}(0,T)}    \| (1+t)^{\frac{3^-}{4}} u \|_{\Ld^{8/3}(0,T; \Ld^\infty(\R^3))}^{p-2} \\
&\lesssim_0  \varphi(\E(0)) \| (1+t)^{p(\gamma-1)    -p\frac{3}{4}^- }\|_{\Ld^{\frac{8}{14-3p}}(0,T)},
\end{aligned}
$$
where we have applied Corollary~\ref{coro-para1} in the last line.
In order to ensure time integrability (and thus a uniform bound independent of $T$), we therefore require that
$$
\left(p \frac{3}{4} - p(\gamma-1)   \right) \frac{8}{14-3p}>1
\Longleftrightarrow
\gamma < \frac{17}{8} - \frac{7}{4p},
$$
which concludes the proof.

\end{proof}

\begin{lem}
\label{lem-source2}
There is $\eps>0$ such that, for all $p \in (3,3+\eps)$, the following holds.  For all $t \in [0,T)$, we have
\begin{equation}
\label{eq:source2-1}
\|  u\cdot \na_x u(t) \|_{\Ld^p(\R^3)} \lesssim_0  \varphi({\E(0)})  \| \Delta_x u (t) \|_{\Ld^p(\R^3)}.
\end{equation}
Consequently, for all $\gamma \geq 0$,
\begin{equation}
\label{eq:source2-2}
\| (1+t)^{\gamma} u\cdot \na_x u \|_{\Ld^p(0,T; \Ld^p(\R^3))} \lesssim_0 \varphi({\E(0)}) \| (1+t)^{\gamma} \Delta_x u \|_{\Ld^p(0,T; \Ld^p(\R^3))}.
\end{equation}

\end{lem}

\begin{proof}
By the Holder inequality, we have
\begin{equation}
\label{eq-hold1}
\|  u\cdot \na_x u \|_{\Ld^p(\R^3)} \leq \|  u\|_{\Ld^{6}(\R^3)} \|  \na_x u \|_{\Ld^{q}(\R^3)}. 
\end{equation}
with $\frac{1}{p} = \frac{1}{6}+\frac{1}{q}$.  
By the Gagliardo-Nirenberg inequality, there are $1>\alpha_1> \alpha_2 >0$ such that
\begin{equation}
\label{eq-gn1}
 \| \na_x  u\|_{\Ld^{q}(\R^3)} \lesssim \| \Delta_x u \|_{\Ld^p(\R^3)}^{\alpha_1} \| u \|_{\Ld^2(\R^3)}^{1-\alpha_1}, \qquad  \|   u\|_{\Ld^{6}(\R^3)} \lesssim \| \Delta_x u \|_{\Ld^p(\R^3)}^{\alpha_2} \| u \|_{\Ld^2(\R^3)}^{1-\alpha_2},
\end{equation}
with
\begin{align*}
\frac{1}{q} = \frac{1}{3}+ \left(\frac{1}{p}-\frac{2}{3}\right)\alpha_1 + \frac{1-\alpha_1}{2}, \qquad
\frac{1}{6} = \left(\frac{1}{p}-\frac{2}{3}\right)\alpha_2 + \frac{1-\alpha_2}{2}.
\end{align*}
We can check that when $p$ is close to $3$, $q$ is close to $6$ while $\alpha_1$ and $\alpha_2$ are close to $4/5$ and $2/5$.
Therefore, we must have 
$$\alpha_1 + \alpha_2 >1.$$
We then write
$$
\|  u\|_{\Ld^{6}(\R^3)} = \|  u\|_{\Ld^{6}(\R^3)}^{\frac{1-\alpha_1}{\alpha_2}}  \|  u\|_{\Ld^{6}(\R^3)}^{1-\frac{1-\alpha_1}{\alpha_2}}.
$$
As by Proposition~\ref{propo:estvns3D}, $\| u \|_{\Ld^\infty(0,T; \Ld^6(\R^3))} \lesssim_0 1$, we have by~\eqref{eq-gn1} the bound
$$
\begin{aligned}
\|  u\|_{\Ld^{6}(\R^3)} &\lesssim_0  \|  u\|_{\Ld^{6}(\R^3)}^{\frac{1-\alpha_1}{\alpha_2}}  \\
&\lesssim_0 \| \Delta_x u \|_{\Ld^p(\R^3)}^{1-\alpha_1} \| u \|_{\Ld^2(\R^3)}^{\frac{(1-\alpha_1)(1-\alpha_2)}{\alpha_2}}.
\end{aligned}
$$
We finally get by~\eqref{eq-hold1},~\eqref{eq-gn1}  and the fact that  $\| u(t) \|_{\Ld^2(\R^3)}^2 \leq 2 \E(0)$,
$$
\begin{aligned}
\|  u\cdot \na_x u \|_{\Ld^p(\R^3)} &\lesssim_0  \| \Delta_x u \|_{\Ld^p(\R^3)} \| u \|_{\Ld^2(\R^3)}^{1-\alpha_1}  \\
&\lesssim_0  \varphi({\E(0)})  \| \Delta_x u \|_{\Ld^p(\R^3)}
\end{aligned}
$$
and the proof of~\eqref{eq:source2-1}
 is complete. 
\end{proof}

\subsection{Weighted maximal parabolic regularity  using higher decay of higher dissipation}

We now turn to the key estimates of the proof.
Namely, we apply maximal regularity in $\Ld^p$ for  some $p>3$ after relying on the higher decay of higher dissipation provided by Lemma~\ref{lem-gaindecay} and Corollary~\ref{coro-Dpmain}.

\begin{lem}
\label{lem-source4}
Let $p>3$.
For all $\gamma\geq 0$, there is $k \in \N$ large enough, so that
\begin{equation}
\int_0^t \D_p (s) (1+s)^{p \gamma-k} \, \dd s \lesssim_0 1 + \| (1+t)^\gamma \Delta u \|_{\Ld^p(0,T; \Ld^p(\R^3))}^{\alpha p}, 
\end{equation}
for some $\alpha \in (0,1)$.
\end{lem}

\begin{proof}
We use the rough bound
\begin{align*}
\int_0^t \D_p (s) (1+s)^{p \gamma-k} \, \dd s 
 &\lesssim \int_0^t  \left(\int_{\R^3 \times \R^3} f(s,x,v) |v|^p  \, \dd v \, \dd x \right)(1+s)^{p\gamma-k} \, \dd s  \\ &\qquad\qquad + \int_0^t  \left( \int_{\R^3} \rho_f(s,x) |u|^p  \,   \dd x \right)(1+s)^{p\gamma-k} \, \dd s \\
  &\lesssim \int_0^t  m_p(s) (1+s)^{p\gamma-k} \, \dd s  +  \int_0^t \| u\|^p_{\Ld^\infty(\R^3)} (1+s)^{p\gamma-k} \, \dd s . 
\end{align*}
The first term above is treated with Lemma~\ref{lem-decaymp}. For the second one, we rely on the same interpolation procedure as in the previous proofs. By the Gagliardo-Nirenberg inequality, we have
$$
\| u\|_{\Ld^\infty(\R^3)} \lesssim  \| \Delta_x u \|_{\Ld^p(\R^3)}^\alpha \| u \|_{\Ld^2(\R^3)}^{1-\alpha},
$$
with 
$$
0= 0+ \left( \frac{1}{p} - \frac{2}{3} \right) \alpha + \frac{1-\alpha}{2} \Longleftrightarrow \alpha = \frac{\frac{1}{2}}{\frac{7}{6}- \frac{1}{p} } = \frac{3p}{7p-6}.
$$ 
We may therefore use the Holder inequality to get
$$
\begin{aligned}
\int_0^t \| u\|^p_{\Ld^\infty(\R^3)} (1+s)^{p\gamma-k} \, \dd s   &\lesssim_0  \int_0^t \| \Delta_x u \|_{\Ld^p(\R^3)}^{\alpha p} (1+s)^{p\gamma-k} \, \dd s \\
&\lesssim_0  \| (1+s)^{ p(1-\alpha) \gamma -k }\|_{\Ld^{\frac{1}{1-\alpha}}(0,T)} \| (1+s)^\gamma \Delta u \|_{\Ld^p(0,T; \Ld^p(\R^3))}^{\alpha p}
\end{aligned}
$$
and take $k$ large enough to ensure uniform integrability in time, so that we get the claimed bound.

\end{proof}

\begin{lem}
\label{lem-source3}
Let $p>3$. For all $\gamma \in \left(0,\frac{27}{8}- \frac{13}{4p}\right)$, there holds
\begin{equation}
\| (1+s)^\gamma |\na_x u| m_p^{1/p} \|_{\Ld^p(0,T; \Ld^p(\R^3))} \lesssim_0 \varphi(\E(0)) \| (1+t)^{\gamma} \Delta_x u \|_{\Ld^p(0,T; \Ld^p(\R^3))}^\alpha,
\end{equation}
for some $\alpha \in (0,1)$.

\end{lem}

\begin{proof}
By Lemma~\ref{lem-decaymp}, we have
$$
(1+t)^{ p \gamma} \||\na_x u| m_p^{1/p}  \|_{\Ld^p(\R^3)}^p \lesssim_0 (1+t)^{ p \left(\gamma - \frac{3}{4}^-\right)}  \|\na_x u  \|_{\Ld^p(\R^3)}^p.
$$
By the Gagliardo-Nirenberg inequality, we can write that
$$
 \|\na_x u  \|_{\Ld^p(\R^3)} \lesssim \| \Delta_x u \|_{\Ld^p(\R^3)}^\alpha \| u \|_{\Ld^2(\R^3)}^{1-\alpha}
$$
with 
$$
\frac{1}{p} = \frac{1}{3} + \left( \frac{1}{p} - \frac{2}{3} \right) \alpha + \frac{1-\alpha}{2}
\Longleftrightarrow 
\alpha = \frac{\frac{5}{6}- \frac{1}{p} }{\frac{7}{6}- \frac{1}{p} } = \frac{5p-6}{7p-6}.
$$ 
Note that since $p>3$, we have $\alpha>3/5$. As a result, by Theorem~\ref{thm-cond}, we obtain
\begin{align*}
(1+t)^{ p \gamma} \||\na_x u| m_p^{1/p}  \|_{\Ld^p(\R^3)}^p  &\lesssim_0 \varphi(\E(0)) (1+t)^{ p \left(\gamma(1-\alpha) + (\alpha -2) \frac{3}{4}^-  \right)} [(1+t)^{\gamma} \| \Delta_x u \|_{\Ld^p(\R^3)}]^{\alpha p}.
\end{align*}
Applying the Holder inequality, we get
\begin{align*}
&\| (1+t)^\gamma |\na_x u| m_p^{1/p} \|_{\Ld^p(0,T; \Ld^p(\R^3))} \|_{\Ld^p(\R^3)}^p \\
  &\qquad \lesssim_0  \varphi(\E(0)) \| (1+t)^{ p \left(\gamma(1-\alpha) + (\alpha -2) \frac{3}{4}^- \alpha \right)} \|_{\Ld^{\frac{1}{1-\alpha}} (0,T)} \| (1+t)^{\gamma} \Delta_x u \|_{\Ld^p(0,T; \Ld^p(\R^3))}^{p \alpha}.
\end{align*}
In order to ensure integrability in time, we thus need to enforce that
$$
p \left( \frac{3}{4}(2-\alpha) - \gamma(1-\alpha)\right)  \frac{1}{1-\alpha} >1 \Longleftrightarrow
\gamma <  \frac{27}{8}- \frac{13}{4p},    
$$
and we obtain the claimed interval for $\gamma$.
\end{proof}

The preceding estimates are considered in order to estimate
$$S= (1+t)^{\gamma/2} (j_f - \rho_f u)  - (1+t)^{\gamma/2} u\cdot \na_x u  + \frac{\gamma}{2} (1+t)^{\gamma/2-1}  u$$
for a sufficiently large value of $\gamma>0$, in $\Ld^p$ for some $p>3$. This is the purpose of the next lemma.
\begin{lem}
\label{lem-sourcefinal}
There is $\eps>0$ such that, for all $p \in (3,3+\eps)$, the following holds. For all  $\gamma \in \left(0,  \frac{17}{8} - \frac{7}{4p} \right)$ we have
\begin{equation}
\label{eq-sourcefinal}
\begin{aligned}
\| S\|_{\Ld^p(0,T; \Ld^p(\R^3))} &\lesssim  \| \rho_f \|_{\Ld^\infty(0,t; \Ld^\infty(\R^3)}  \| (1+s)^\gamma \pa_s u \|_{\Ld^p(0,t; \Ld^p(\R^3)} \\
&\quad+ \psi\left(\|u_0\|_{\H^{1/2}(\R^3)}+M_\alpha f_0+N_q(f_0)+\E(0) + 1 \right) \\
&\times\Bigg[
1 + \| (1+t)^\gamma \Delta u \|_{\Ld^p(0,T; \Ld^p(\R^3))}^{\alpha p} + \varphi(\E(0)) \| (1+t)^\gamma \Delta u \|_{\Ld^p(0,T; \Ld^p(\R^3))}^{p} \Bigg]^{1/p},
\end{aligned}
\end{equation}
where $\alpha \in (0,1)$ appears in Lemma~\ref{lem-source4}.
\end{lem}

\begin{proof}
For $\gamma<  \frac{17}{8} - \frac{7}{4p}$, the contribution of the terms ${\gamma} (1+t)^{\gamma-1}  u$ and $(1+t)^{\gamma} u\cdot \na_x u$ are treated thanks to  Lemma~\ref{lem-source1} and  Lemma~\ref{lem-source2}, 

The contribution of the Brinkman force requires the use of the higher decay of higher dissipation as provided by Lemma~\ref{lem-Dp} and Lemma~\ref{lem-gaindecay}.
By Corollary~\ref{coro-Dpmain} (more precisely by~\eqref{eq-keyDp}),
we first  bound
\begin{align*}
\| (1+t)^\gamma (j_f - \rho_f u)\|_{\Ld^p(0,T; \Ld^p(\R^3))}^p 
   & \lesssim \| \rho_f \|_{\Ld^\infty(0,T;\Ld^\infty(\R^3))}^{{p-1}}  \int_0^T  \D_p(t)  (1+t)^{p\gamma}  \, \dd t \\
& \lesssim  \| \rho_f \|_{\Ld^\infty(0,T;\Ld^\infty(\R^3))}^{{p-1}} \Bigg( \int_0^t \D_p (s) (1+s)^{p \gamma-k} \, \dd s \\
&  \quad + \| \rho_f \|_{\Ld^\infty(0,t; \Ld^\infty(\R^3))}  \| (1+s)^\gamma \pa_s u \|_{\Ld^p(0,t; \Ld^p(\R^3)}^p \\
&\quad \quad+ \| (1+s)^\gamma |\na_x u| m_p^{1/p} \|_{\Ld^p(0,t; \Ld^p(\R^3)}^p  + \int_{\R^3 \times \R^3} f_0 |v-u_0|^p \, \dd v \, \dd x\Bigg).
\end{align*}

Recall the meaning of the notation $\lesssim_0$ from Definition~\ref{nota:lesssim}.
By Lemma~\ref{lem-source4} and Lemma~\ref{lem-source3}, we infer that for $\gamma<\frac{27}{8}- \frac{13}{4p}$ and $k$ large enough,  
\begin{align*}
&\| (1+t)^\gamma (j_f - \rho_f u)\|_{\Ld^p(0,T; \Ld^p(\R^3))}^p   \\
&\lesssim \| \rho_f \|_{\Ld^\infty(0,t; \Ld^\infty(\R^3)}^p  \| (1+s)^\gamma \pa_s u \|_{\Ld^p(0,t; \Ld^p(\R^3)}^p \\
&\quad+ \psi\left(\|u_0\|_{\H^{1/2}(\R^3)}+M_\alpha f_0+N_q(f_0)+\E(0) + 1 \right) \Bigg[
1 + \| (1+t)^\gamma \Delta u \|_{\Ld^p(0,T; \Ld^p(\R^3))}^{\alpha p} \\
&\quad \qquad + \varphi(\E(0)) \| (1+t)^\gamma \Delta u \|_{\Ld^p(0,T; \Ld^p(\R^3))}^{p} \Bigg],
\end{align*}
with $\alpha \in (0,1)$. This concludes the proof of the lemma.

\end{proof}

We are finally in position to apply Theorem~\ref{thm-maximal} for some $p=q>3$.
\begin{lem}
\label{lem-max1-3}
There is $\eps>0$ such that, for all $p \in (3,3+\eps)$, the following holds. For all  $\gamma \in \left(0,  \frac{17}{8} - \frac{7}{4p} \right)$ we have
\begin{equation}
\| (1+t)^{\gamma} \pa_t u \|_{\Ld^p(0,T; \Ld^p(\R^3))} + \| (1+t)^{\gamma} \Delta u \|_{\Ld^p(0,T; \Ld^p(\R^3))} \lesssim_0 1.
\end{equation}
\end{lem}

\begin{proof}
We apply Theorem~\ref{thm-maximal}  together with Lemma~\ref{lem-sourcefinal} and the Young inequality,
which yields
\begin{align*}
&\| (1+t)^{\gamma} \pa_t u \|_{\Ld^p(0,T; \Ld^p(\R^3))} + \| (1+t)^{\gamma} \Delta u \|_{\Ld^p(0,T; \Ld^p(\R^3))}  \\
&\leq \psi_1\left(\|u_0\|_{\H^{1/2}(\R^3)}+M_\alpha f_0+N_q(f_0)+\E(0) + 1 \right) + \| u_0\|_{B^{s_p,p}_p(\R^3)} \\
&\quad+ \Bigg[\psi_2\left(\|u_0\|_{\H^{1/2}(\R^3)}+M_\alpha f_0+N_q(f_0)+\E(0) + 1 \right)\varphi(\E(0))+ \C \| \rho_f \|_{\Ld^\infty(0,T; \Ld^\infty(\R^3))} + \frac{1}{2}\Bigg]  \\
&\qquad \quad  \times \Bigg(\| (1+t)^{\gamma} \pa_t u \|_{\Ld^p(0,T; \Ld^p(\R^3))} + \| (1+t)^{\gamma} \Delta u \|_{\Ld^p(0,T; \Ld^p(\R^3))} \Bigg).
\end{align*}
Recall that by Lemma~\ref{lem:rho}, we have $\| \rho_f \|_{\Ld^\infty(0,t; \Ld^\infty(\R^3))} \lesssim \| f_0 \|_{\Ld^1_v(\R^3; \Ld^\infty_x(\R^3))} $.
We may choose $\Psi$ and $\delta_0$ in the assumptions of Theorem~\ref{thm} such that by~\eqref{eq:introsmall}, both $\E(0)$ and   $\| f_0 \|_{\Ld^1_v(\R^3; \Ld^\infty_x(\R^3))}$ are small enough so that 
$$
\Bigg[\psi_2\left(\|u_0\|_{\H^{1/2}(\R^3)}+M_\alpha f_0+N_q(f_0)+\E(0) + 1 \right)\varphi(\E(0))+ \C \| \rho_f \|_{\Ld^\infty(0,T; \Ld^\infty(\R^3))} + \frac{1}{2}\Bigg]<\frac{2}{3}.$$
We can thus absorb all terms of the rhs involving $\| (1+t)^{\gamma} \pa_t u \|_{\Ld^p(0,T; \Ld^p(\R^3))} + \| (1+t)^{\gamma} \Delta u \|_{\Ld^p(0,T; \Ld^p(\R^3))} $ by the lhs. The proof of the lemma is therefore complete.
\end{proof}

 As we can take $\gamma>1/p'$, 
arguing again by interpolation, we obtain the desired estimate for $\| \nabla_x u\|_{\Ld^\infty(\R^3)}$. 

\begin{coro}
\label{coro-para1-3}
There holds
$$
\int_0^T \| \na u \|_{\Ld^\infty(\R^3)} \, \dd s \lesssim_0 \varphi(\E(0)).
$$
\end{coro}

\begin{proof} By the Gagliardo-Nirenberg inequality, for $p >3 $,
$$
\| \na u \|_{\Ld^\infty(\R^3)} \lesssim \| \Delta u \|_{\Ld^p(\R^3)}^\beta \| u \|_{\Ld^2(\R^3)}^{1-\beta}
$$
with
$$
0 = \frac{1}{3} + \left(\frac{1}{p}-\frac{2}{3} \right) \beta +\frac{1-\beta}{2} \Longleftrightarrow \beta = \frac{5p}{7p-6}.
$$
We choose $p$ as in the statement of Lemma~\ref{lem-max1-3}. Note that for $p$ close to $3$, $\beta$ is close to $1$.
By the Holder inequality, we thus obtain
\begin{align*}
\int_0^T \| \na u \|_{\Ld^\infty(\R^3)} \, \dd s &\lesssim \int_0^T  \| \Delta u \|_{\Ld^p(\R^3)}^\beta \| u \|_{\Ld^2(\R^3)}^{1-\beta} \, \dd s \\
&\lesssim  \int_0^T  \| \Delta u \|_{\Ld^p(\R^3)}^\beta \| u \|_{\Ld^2(\R^3)}^{1-\beta} \, \dd s \\
&\lesssim  \int_0^T  \left[(1+s)^\gamma\| \Delta u \|_{\Ld^p(\R^3)}\right]^\beta  (1+s)^{-\gamma \beta} \| u \|_{\Ld^2(\R^3)}^{1-\beta} \, \dd s \\
&\lesssim  \left(\int_0^T (1+t)^{-\gamma \frac{\beta p}{p-\beta} } \| u \|_{\Ld^2(\R^3)}^{(1-\beta) \frac{p}{p-\beta}} \, \dd s \right)^{ \frac{p-\beta}{p}}  \| (1+t)^{\gamma} \Delta u \|_{\Ld^p(0,T; \Ld^p(\R^3))}^\beta.
\end{align*}
As we can choose $p$ as close to $3$ as necessary (by choosing $p_0$ appropriately in the assumptions of Theorem~\ref{thm}) and take  $\gamma$ close to  $\frac{17}{8} - \frac{7}{4p}$ thanks to Lemma~\ref{lem-max1-3}, we can ensure
$$
\gamma \frac{\beta p}{p-\beta}>1,
$$
and we may finally combine it with the energy bound $ \| u \|_{\Ld^2(\R^3)}^2 \lesssim \E(0)$ to get the claimed result.

\end{proof}

\subsection{End of the bootstrap}

We are in position to conclude.

\begin{lem}
We have
\begin{equation}
\label{eq:H12-petit}
\int_0^{t^\star} \| j_f - \rho_f u \|_{\dot{\H}^{1/2}(\R^3)}^2 \, \dd t  \lesssim_0 \E(0).
\end{equation}
\end{lem}

\begin{proof}By the Sobolev embedding and the Holder inequality,
$$
\begin{aligned}
 \| j_f - \rho_f u \|_{\dot{\H}^{1/2}(\R^3)}^2&\lesssim \| j_f - \rho_f u \|_{\Ld^{3/2}(\R^3)}^2 \\
 &\lesssim \D_2 \| \rho_f \|_{\Ld^3(\R^3)}.
 \end{aligned}
$$
Using the energy inequality and the fact that 
$$
\sup_{[0,t^\star]}  \| \rho_f \|_{\Ld^3(\R^3)} \leq \sup_{[0,t^\star]}  \| \rho_f \|_{\Ld^1(\R^3)}^{1/3}   \| \rho_f \|_{\Ld^{\infty}(\R^3)}^{2/3}    \lesssim_0  1,
$$
we finally obtain
$$
\int_0^{t^\star} \| j_f - \rho_f u \|_{\dot{\H}^{1/2}(\R^3)}^2 \, \dd t \lesssim_0 \int_0^{t^\star} \D_2(t) \, \dd t  \lesssim_0 \E(0),
$$
concluding the proof of~\eqref{eq:H12-petit}. 

\end{proof}
We deduce that choosing $\Psi$ appropriately, \eqref{eq:introsmall} enforces that
\begin{align}
\|u_0\|_{\dot \H^{1/2}(\R^3)}^2  + \C_\star \int_0^{t^\star} \|F(s)\|_{\dot \H^{-1/2}(\R^3)}^2\,\dd s < \frac{1}{\C_\star^2}.
\end{align}
Recalling~\eqref{eq:H12-general}, this means that there exist strong existence times that are strictly larger than $t^\star$.
On the other hand, choosing $\Psi$ appropriately for~\eqref{eq:introsmall}, Corollary~\ref{coro-para1-3} entails that 
$$
\int_0^{t^\star} \| \na u \|_{\Ld^\infty(\R^3)} \, \dd s \leq  \frac{\delta_0}{2},
$$
where $\delta_0$ is the parameter of Lemma~\ref{charac}.
Owing to~\eqref{nablaulocal} in Proposition~\ref{propo:local},  we can thus find a strong existence time $t_0 > t^\star$ such that 
$$
\int_0^{t_0} \| \na u \|_{\Ld^\infty(\R^3)} \, \dd s <  {\delta_0}.
$$
This is a contradiction with the definition of $t^\star$.
We deduce that we must have $t^\star= +\infty$ and the proof of Theorem~\ref{thm} is complete.

\bigskip 
\noindent {\bf Acknowledgements.} Partial support by the grant ANR-19-CE40-0004 is acknowledged. I warmly  thank Ayman Moussa and Iv\'an Moyano for numerous discussions about the Vlasov-Navier-Stokes system over the past years and for our joint paper \cite{HKMM} without which this work wouldn't have been possible.

\bibliographystyle{abbrv}
\bibliography{vns}
\end{document}